\documentclass[a4paper]{article}
\usepackage{amsfonts}
\usepackage{amsmath}
\usepackage{graphicx}

\providecommand{\U}[1]{\protect\rule{.1in}{.1in}}
\newtheorem{theorem}{Theorem}

\newtheorem{corollary}[theorem]{Corollary}

\newtheorem{definition}[theorem]{Definition}

\newtheorem{lemma}[theorem]{Lemma}

\newtheorem{proposition}[theorem]{Proposition}
\newtheorem{remark}[theorem]{Remark}

\newenvironment{proof}[1][Proof]{\noindent\textbf{#1.} }{\ \rule{0.5em}{0.5em}}
\input{tcilatex}
\begin{document}

\title{Polyharmonic Hardy Spaces on the Klein-Dirac Quadric with Application
to Polyharmonic Interpolation and Cubature Formulas }
\author{Ognyan Kounchev and Hermann Render}
\maketitle

\begin{abstract}
In the present paper we introduce a new concept of Hardy type space
naturally defined on the Klein-Dirac quadric. We study different properties
of the functions belonging to these spaces, in particular boundary value
problems. We apply these new spaces to polyharmonic interpolation and to
interpolatory cubature formulas.
\end{abstract}

\section{ Introduction}

In one-dimensional mathematical analysis, Interpolation Theory and
Quadrature formulas are intimately related, cf. \cite{krylov}, \cite{davis}.
This relation causes a similarity between the approaches for estimation of
the remainders of Interpolation and Quadrature. One approach for estimation
of the error in Interpolation theory is related to Lagrange formula and uses
higher dierivatives of the interpolated function (cf. \cite{krylov}, chapter 
$3,$ Theorem $4,$ and \cite{davis}, Theorem $3.1.1$). The second approach
uses analyticity of the interpolated function and Hermite formula (cf. \cite%
{krylov}, chapter $3,$ Theorem $5,$ and \cite{davis}, Theorema $3.6.1$). In
a similar way, already A. Markov has estimated the error of a quadrature
formula for differentiable functions in $C^{N}\left( I\right) $ defined on
the interval $I$ by means of its $N-$th derivative (cf. \cite{krylov},
chapter $7.1$, and Davis \cite{davis}, p. $344$). The second approach
estimates the error of a quadrature formula for certain classes of functions 
$f$ which are \emph{analytic} on some open set $D$ in $\mathbb{C}$
containing the interval $I,$ (cf. \cite{davisRabinowitz}, chapter $4.6,$ see
also \cite{krylov}, chapter $12.2$). However, in both Interpolation and
Quadrature, the first approach is usually not very practical beyond
derivatives of order five, see \cite{stroud}.

Interpolation and Approximation of integrals in the multivariate case is a
much more difficult task. In Numerical Analysis, instead of quadrature
formula the notion of cubature formula is often used, see \cite{sobolev}, 
\cite{stroudBook}, \cite{sobolev2} and the recent survey \cite{cools}. In
contrast to the univariate case there is no satisfactory error analysis
available in the multivariate case, cf. \cite{stroudBook}, \cite{bakhvalov},
and part $4$ in the last Russian edition of the classical monograph \cite%
{krylov}. Let us mention that the area of quadrature domains which has
received a lot of interest recently presents an interesting multidimensional
alternative and we refer to \cite{ebenfelt}.

The present research continues the study of estimates of \emph{polyharmonic
interpolation }and \emph{polyharmonic interpolatory cubature formulas}
initiated in \cite{kounchevRenderHardyAnnulus}; in the present paper we
consider interpolation and cubature formulas in the ball in $\mathbb{R}^{d},$
while in \cite{kounchevRenderHardyAnnulus} the case of an annular region was
considered. 

\subsection{Gauss-Almansi formula}

In \cite{haussmannKounchev} \emph{polyharmonic interpolation} has been
considered for functions defined in the ball in $\mathbb{R}^{d}.$ In the
same spirit, in \cite{kounchevRenderArkiv} and \cite{kounchevRenderArxiv} we
have introduced a new multivariate cubature formulae $C_{N}\left( f\right) $
in the ball depending on a parameter $N\in \mathbb{N}$ which approximates
the integral 
\begin{equation}
\dint\limits_{B_{R}}f\left( x\right) d\mu \left( x\right) 
\label{eqintegral}
\end{equation}%
for continuous functions $f:B_{R}\rightarrow \mathbb{C}$ defined on the 
\emph{ball} 
\begin{equation}
B_{R}=\left\{ x\in \mathbb{R}^{d}:\left\vert x\right\vert <R\right\} ,
\label{Aab}
\end{equation}%
where $\left\vert x\right\vert $ denotes the euclidean norm of $x=\left(
x_{1},...,x_{d}\right) \in \mathbb{R}^{d}.$

The exact definition of the \emph{polyharmonic interpolation} formula and of
the \emph{polyharmonic cubature} formula $C_{N}\left( f\right) $ will be
explained in Section \ref{Scubature}. A major purpose of the present paper
is to provide an error analysis for a class of functions on the ball $B_{R}$
which exhibit a certain type of analytical behavior.

Let us introduce the necessary notions and notations. Let 
\begin{equation*}
\mathbb{S}^{d-1}:=\left\{ x\in \mathbb{R}^{d}:\left\vert x\right\vert
=1\right\} 
\end{equation*}%
be the unit sphere endowed with the rotation invariant measure $d\theta $.
We shall write $x\in \mathbb{R}^{d}$ in spherical coordinates $x=r\theta $
with $\theta \in \mathbb{S}^{d-1}.$ Let $\mathcal{H}_{k}\left( \mathbb{R}%
^{d}\right) $ be the set of all harmonic homogeneous complex-valued
polynomials of degree $k.$ Then $f\in \mathcal{H}_{k}\left( \mathbb{R}%
^{d}\right) $ is called a \emph{solid harmonic} and the restriction of $f$
to $\mathbb{S}^{d-1}$ a \emph{spherical harmonic} of degree $k$ and we set 
\begin{equation}
a_{k}:=\dim \mathcal{H}_{k}\left( \mathbb{R}^{d}\right) ,  \label{eqdim}
\end{equation}%
see \cite{steinWeiss}, \cite{seeley}, \cite{askey}, \cite{okbook} for
details. Throughout the paper we shall assume that the set of functions 
\begin{equation}
Y_{k,\ell }:\mathbb{R}^{d}\rightarrow \mathbb{C},\ \text{for }\ell
=1,...,a_{k},  \label{Ykl}
\end{equation}%
is an \emph{orthonormal basis} of $\mathcal{H}_{k}\left( \mathbb{R}%
^{d}\right) $ with respect to the scalar product 
\begin{equation*}
\left\langle f,g\right\rangle _{\mathbb{S}^{d-1}}:=\int_{\mathbb{S}%
^{d-1}}f\left( \theta \right) \overline{g\left( \theta \right) }d\theta .
\end{equation*}%
Recall that due to the homogeneity of $Y_{k,\ell }\left( x\right) $ we have
the identity $Y_{k,\ell }\left( x\right) =r^{k}Y_{k\ell }\left( \theta
\right) $ for $x=r\theta .$

Our polyharmonic Interpolation and polyharmonic Cubature $C_{N}\left(
f\right) $ approximating the integral (\ref{eqintegral}) are based on the 
\emph{Laplace--Fourier series }of the continuous function $%
f:B_{R}\rightarrow \mathbb{C}$, defined by the formal expansion 
\begin{equation}
f\left( r\theta \right) =\sum_{k=0}^{\infty }\sum_{\ell =1}^{a_{k}}f_{k,\ell
}\left( r\right) Y_{k,\ell }\left( \theta \right)  \label{frtheta}
\end{equation}
where the \emph{Laplace--Fourier coefficient }$f_{k,\ell }\left( r\right) $%
\emph{\ }is defined by 
\begin{equation}
f_{k,\ell }\left( r\right) =\int_{\mathbb{S}^{d-1}}f\left( r\theta \right)
Y_{k,\ell }\left( \theta \right) d\theta  \label{eqfourier}
\end{equation}
for any positive real number $r$ with $r<R$ and $a_{k}$ is defined in (\ref%
{eqdim}). There is a strong interplay between algebraic and analytic
properties of the function $f$ and those of the Laplace-Fourier coefficients 
$f_{k,\ell }$. For example, if $f\left( x\right) $ is a polynomial in the
variable $x=\left( x_{1},...,x_{d}\right) $ then the Laplace-Fourier
coefficient $f_{k,\ell }$ is of the form $f_{k,\ell }\left( r\right)
=r^{k}p_{k,\ell }\left( r^{2}\right) $ where $p_{k,\ell }$ is a univariate
polynomial, see e.g. in \cite{steinWeiss} or \cite{sobolev}. Hence, the 
\emph{Laplace-Fourier} series (\ref{frtheta}) of a polynomial $f\left(
x\right) $ is equal to 
\begin{equation}
f\left( x\right) =\sum_{k=0}^{\deg f}\sum_{\ell =1}^{a_{k}}p_{k,\ell
}(\left| x\right| ^{2})Y_{k,\ell }\left( x\right) =\sum_{k=0}^{\deg
f}\sum_{\ell =1}^{a_{k}}\left| x\right| ^{k}p_{k,\ell }(\left| x\right|
^{2})Y_{k,\ell }\left( \theta \right)  \label{gauss}
\end{equation}
where $\deg f$ is the total degree of $f$ and $p_{k,\ell }$ is a univariate
polynomial of degree $\leq \deg f-k.$ This representation is often called
the \emph{Gauss representation.} A similar formula is valid for a much
larger class of functions. Let us recall that a function $f:G\rightarrow 
\mathbb{C}$ defined on an open set $G$ in $\mathbb{R}^{d}$ is called \emph{%
polyharmonic of order} $N$ if $f$ is $2N$ times continuously differentiable
and 
\begin{equation}
\Delta ^{N}u\left( x\right) =0  \label{dlNux=0}
\end{equation}
for all $x\in G$ where $\Delta =\frac{\partial ^{2}}{\partial x_{1}^{2}}+...+%
\frac{\partial ^{2}}{\partial x_{d}^{2}}$ is the Laplace operator and $%
\Delta ^{N}$ the $N$-th iterate of $\Delta .$ The theorem of \emph{Almansi}
states that for a polyharmonic function $f$ of order $N$ defined on the ball 
$B_{R}=\left\{ x\in \mathbb{R}^{d}:\left| x\right| <R\right\} $ there exist
univariate polynomials $p_{k,\ell }\left( r\right) $ of degree $\leq N-1$
such that 
\begin{equation}
f\left( x\right) =\sum_{k=0}^{\infty }\sum_{\ell =1}^{a_{k}}p_{k,\ell
}(\left| x\right| ^{2})Y_{k,\ell }\left( x\right) =\sum_{k=0}^{\infty
}\sum_{\ell =1}^{a_{k}}\left| x\right| ^{k}p_{k,\ell }(\left| x\right|
^{2})Y_{k,\ell }\left( \theta \right)  \label{eqAlmansi}
\end{equation}
where convergence of the sum is uniform on compact subsets of $B_{R},$ see
e.g. \cite{sobolev}, \cite{avanissian}, \cite{aron}.

To end this technical introduction, let us remind some estimates for
spherical harmonics which we will need below:\ 

\begin{enumerate}
\item For every multiindex $\alpha $ and for every integer $k\geq 1$ holds 
\begin{equation}
\left\vert D_{\theta }^{\alpha }Y_{k,\ell }\left( \theta \right) \right\vert
\leq Ck^{\left\vert \alpha \right\vert +\frac{d-2}{2}}\qquad \text{for }%
\theta \in \mathbb{S}^{d-1},  \label{SphericalEstimate}
\end{equation}%
see \cite{seeley}, p. $120.$ Here $D_{\theta }^{\alpha }$ is the multi-index
notation for the derivative with respect to $\theta \in \mathbb{S}^{d-1}.$ 

\item A function $f\left( \theta \right) $ defined on $\mathbb{S}^{d-1}$ is
real analytic if its Laplace-Fourier expansion $f\left( \theta \right)
=\dsum_{k=0}^{\infty }Y_{k}\left( \theta \right) $ (where we have put $%
Y_{k}\left( \theta \right) =\dsum_{\ell =1}^{a_{k}}f_{k,\ell }Y_{k,\ell
}\left( \theta \right) $ ) satisfies 
\begin{equation*}
\left\Vert Y_{k}\left( \theta \right) \right\Vert _{L_{2}\left( \mathbb{S}%
^{d-1}\right) }<Ce^{-\eta k}\qquad \text{for }k\geq 0,
\end{equation*}%
for some constants $C,$ $\eta >0;$ see \cite{sobolev}.
\end{enumerate}

\subsection{Complexification of the ball in $\mathbb{R}^{d},$ related to the
ball of the Klein-Dirac quadric}

We want to study analytical extensions of functions $f$ defined on the ball
using the Laplace-Fourier series (\ref{eqAlmansi}). Our strategy is to
require minimal assumptions on the functions $f;$ thus instead of the
standard approach where one works with functions $f$ which are \emph{a priori%
} analytically extendible to a fixed domain $U$ in the complex space $%
\mathbb{C}^{d}$ (as in \cite{aron}) we shall require only that we can extend
the function $x=r\theta \longmapsto f\left( r\theta \right) $ to an analytic
function $z\theta \longmapsto f\left( z\theta \right) ,$ so we only
complexify the radial variable $r$ to a complex variable $z.$ Henceforth we
will use the following \textbf{terminological convention}:\ If the function $%
f\left( r\theta \right) $ possesses an analytic extension with respect to $r$
we call the extended function $f\left( z\theta \right) $ "$r-$\textbf{%
analytic complexification}" or \textquotedblright $r-$\textbf{analytic
continuation}\textquotedblright . Thus, for the $r-$complexification of the
function $f$ one should expect from equation (\ref{eqfourier}) that the
Laplace-Fourier coefficient $f_{k,\ell }\left( r\right) $ extends to an
analytic function of one variable. Hence, we consider the analytically
continued functions on domains in the set 
\begin{equation*}
\mathbb{C}\times \mathbb{S}^{d-1}.
\end{equation*}

We obtain the following important proposition.

\begin{proposition}
\label{PPolyhBasis} The set of functions 
\begin{equation}
B_{k,\ell ;N}:=\left\{ 
\begin{array}{c}
b_{k,\ell ;j}\left( r,\theta \right) =r^{k+2j}Y_{k,\ell }\left( \theta
\right) : \\ 
\quad k\geq 0,\ \ell =1,2,...,a_{k},\ j=0,1,...,N-1%
\end{array}%
\right\}   \label{PolyhBasis}
\end{equation}%
is \textbf{a basis} for the multivariate polynomials which are polyharmonic
of order $N.$ 
\end{proposition}

The proof follows by representation (\ref{gauss}).

\begin{remark}
We will see later that the basis  $B_{k,\ell ;N}$ is a natural
generalization of the basis $\left\{ r^{j}\right\} _{j=0}^{N-1}$ for the
polynomials in the one-dimensional case.
\end{remark}

The main approach in the present paper is to consider the $r-$%
complexification of the functions $f$ defined on the ball $B_{R}\subset 
\mathbb{R}^{d}$  in the form $f\left( z\theta \right) .$ In particular, for
polyharmonic functions $f$ we can provide a formal expression for the $r-$%
complexification, using representation (\ref{eqAlmansi}),  by the following
formula: 
\begin{equation}
f\left( z\theta \right) =\sum_{k=0}^{\infty }\sum_{\ell =1}^{a_{k}}p_{k,\ell
}(z^{2})Y_{k,\ell }\left( z\theta \right) =\sum_{k=0}^{\infty }\sum_{\ell
=1}^{a_{k}}z^{k}p_{k,\ell }(z^{2})Y_{k,\ell }\left( \theta \right) 
\label{eqAlmansiCompexified}
\end{equation}%
The question of convergence will be addressed in the course of the paper.
From this formula we make a \textbf{crucial observation}: the
complexification $f\left( z\theta \right) $ depends only on the values $%
z\theta $ and not on the coordinates of the pair $\left( z,\theta \right) .$
Indeed, this is due to the equality $z^{2}=\left( z\theta ,z\theta \right) $
where for the vectors $w,u\in \mathbb{C}^{d}$ we have the \emph{non-Hermitian%
} product by putting $\left( u,w\right) :=\dsum_{j=1}^{d}u_{j}w_{j}.$ Thus
the function $f\left( z\theta \right) $ is defined on the space 
\begin{equation*}
\mathbb{C}\times \mathbb{S}^{d-1}/\mathbb{Z}_{2}=\left\{ z\theta :z\in 
\mathbb{C},\ \theta \in \mathbb{S}^{d-1}\right\} ,
\end{equation*}%
where the factor in $\mathbb{Z}_{2}$ means identification of the points $%
\left( z,\theta \right) $ and $\left( -z,-\theta \right) $ in $\mathbb{C}%
\times \mathbb{S}^{d-1}.$ But the set $\mathbb{C}\times \mathbb{S}^{d-1}/%
\mathbb{Z}_{2}\subset \mathbb{C}^{d}$ is one of the possible representations
of the famous \emph{Klein-Dirac quadric}.

\begin{definition}
We define the \textbf{Klein-Dirac} quadric by putting 
\begin{equation}
\func{KDQ}:=\mathbb{C}\times\mathbb{S}^{d-1}/\mathbb{Z}_{2}.  \label{KDQ}
\end{equation}
\end{definition}

We shall see that the $r-$analytic continuation of the solutions of the
polyharmonic equation are in fact \textquotedblright the analytic
functions\textquotedblright\ naturally defined on the Klein-Dirac quadric $%
\func{KDQ}.$ The main interest of the present paper is devoted to the
Function theory on the \textbf{complexified ball } $\mathcal{B}_{R}$ in $%
\func{KDQ}$ defined by 
\begin{equation}
\mathcal{B}_{R}:=\left\{ z\theta :\ \left\vert z\right\vert <R,\ \theta \in 
\mathbb{S}^{d-1}\right\} =\mathbb{D}_{R}\times \mathbb{S}^{d-1}/\mathbb{Z}%
_{2},  \label{BallComplexified}
\end{equation}%
where $\mathbb{D}_{R}$ is the open disc of radius $R$ in $\mathbb{C}$, i.e. $%
\mathbb{D}_{R}:=\left\{ z\in \mathbb{C}:\left\vert z\right\vert <R\right\} ;$
for $R=1$ as usually one puts 
\begin{equation*}
\mathbb{D}_{1}=\mathbb{D}\text{ and }\mathcal{B}=\mathcal{B}_{1}.
\end{equation*}

Some enlightening comments about the \textbf{Klein-Dirac quadric} $\func{KDQ}
$ defined in (\ref{KDQ}) are in order. 

\begin{remark}
This quadric has been originally introduced in a special case by Felix Klein
in his Erlangen program in $1870$, where he put forth his correspondence
between the lines in complex projective 3-space and a general quadric in
projective 5-space. The physical relevance of the quadric and the relation
to the conformal motions of compactified Minkowski space-time had been
exploited by Paul Dirac in $1936$ \cite{dirac}. The Klein-Dirac quadric
plays an important role also in Twistor theory, where it is related to the
complexified compactified Minkowski space \cite{penrose}. Apparently, the
term \textquotedblright Klein-Dirac quadric\textquotedblright\ for arbitrary
dimension $d,$ has been coined by the theoretical physicist I. Todorov, cf.
e.g. \cite{nikolovTodorov}, \cite{nikolovTodorov2}. In these references
important aspects of the Function theory on the ball in $\func{KDQ}$ have
been considered in the context of Conformal Quantum Field Theory (CFT). In
the context of CFT Laurent expansions appear in a natural way as the field
functions in the higher dimensional conformal vertex algebras (using a
complex variable parametrization of compactified Minkowski space); see in
particular formula (4.43) in \cite{nikolovTodorov3}, as well as the
references \cite{nikolovTodorov}, \cite{nikolovTodorov2}.
\end{remark}

Our main novelty will be a multivariate generalization of the classical
Hardy space $H^{2}\left( \mathbb{D}\right) $ called \textbf{(polyharmonic)
Hardy space on the ball }$\mathcal{B}_{R}$ to be introduced in Definition %
\ref{DHardyball}. We will denote this space by $H^{2}\left( \mathcal{B}%
_{R}\right) .$ For simplicity sake we will restrict ourselves to considering
the space $H^{2}\left( \mathcal{B}\right) .$

Running ahead of the events, let us say that the name \emph{polyharmonic}
comes from the fact that $H^{2}\left( \mathcal{B}_{R}\right) $ may be
obtained as a limit of the complexifications of the polyharmonic functions
in the ball (\ref{eqAlmansiCompexified}): we take the closure of all finite
sums of the type 
\begin{equation}
u\left( z,\theta \right) =\sum_{k=0}^{\infty }\sum_{\ell
=1}^{a_{k}}u_{k,\ell }\left( z^{2}\right) z^{k}Y_{k,\ell }\left( \theta
\right) ,  \label{uzthetaAlmansiBall}
\end{equation}%
where $u_{k,\ell }\left( \cdot \right) $ are algebraic polynomials of degree 
$\leq N-1;$ such functions $u$ satisfy $\Delta ^{N}u\left( x\right) =0.$ In
view of the Polyharmonic Paradigm announced in \cite{okbook}, the space $%
H^{2}\left( \mathcal{B}_{R}\right) $ generalizes the classical Hardy spaces
which are obtained as limits of algebraic polynomials, where the degree of a
polynomial is replaced by the \textbf{degree of polyharmonicity}.

The Hardy space $H^{2}\left( \mathcal{B}_{R}\right) $ will be a Hilbert
space and we will provide a \emph{Cauchy type kernel}, which is the analog
and a generalization to the Hua-Aronszajn kernel in the ball (cf. \cite{aron}%
, p. 125, Corollary $1.1$). Let us note that the last is a multidimensional
generalization of the classical Cauchy kernel $\frac{1}{z-a}$ from Complex
Analysis (more about Cauchy kernels see in \cite{kerzmanStein}, \cite%
{krantzHarmonic}).

\begin{remark}
The reader familiar with the Cartan classification of classical domains, may
remark that the boundary $\partial \mathcal{B}_{R}$ is very close to the 
\textbf{Shilov boundary} 
\begin{equation*}
\left\{ \left( e^{i\varphi }r\theta \right) :\varphi \in \left[ 0,\pi \right]
,\ 0<r<R,\ \theta \in \mathbb{S}^{d-1}\right\} 
\end{equation*}%
of the so-called \textbf{Cartan classical domai}n $\mathcal{R}_{IV}$ (called
also \textquotedblright Lie-ball\textquotedblright )\ equal to 
\begin{align*}
\widehat{B}_{1}& :=\left\{ \xi +i\eta \in \mathbb{C}^{d}:\xi ,\eta \in 
\mathbb{R}^{d},\ q\left( \xi +i\eta \right) <1\right\} \qquad \text{where }
\\
q\left( \xi +i\eta \right) & =\sqrt{\left\vert \xi \right\vert
^{2}+\left\vert \eta \right\vert ^{2}+2\sqrt{\left\vert \xi \right\vert
^{2}\left\vert \eta \right\vert ^{2}-\left\langle \xi ,\eta \right\rangle
^{2}}}.
\end{align*}%
This has been considered from the point of view of several complex variables
in the monograph of Hua \cite{Hua}, and for the study of the polyharmonic
functions of infinite order in the monographs \cite{aron} (see in particular
p. $59$ and $126$) and \cite{avanissian}. The Hardy spaces defined in
Definition \ref{DHardyball} can be identified with the Hardy space of
holomorphic functions on the Lie ball in $\mathbb{C}^{d},$ cf. \cite{Shai03}%
; this correspondence will be given a thorough consideration in \cite%
{kounchevRenderBook}.
\end{remark}

\begin{remark}
Let us define the annulus in the Klein-Dirac quadric $\func{KDQ}$ as the set
\begin{equation*}
\widetilde{A}_{a,b}:=\left\{ z\theta \in \func{KDQ}:a<\left\vert
z\right\vert <b,\ \theta \in \mathbb{S}^{d-1}\right\} .
\end{equation*}%
The Function theory on $\widetilde{A}_{a,b}$ would help to relate the
present results to previous obtained by us. We have seen in \cite%
{kounchevRenderHardyAnnulus} that an interesting, consistent and fruitful
Function theory is available only on the set 
\begin{equation*}
\mathcal{A}_{a,b}=\left\{ \left( z,\theta \right) :a<\left\vert z\right\vert
<b,\ \theta \in \mathbb{S}^{d-1}\right\} 
\end{equation*}%
which is a subset of $\mathbb{C}\times \mathbb{S}^{d-1}.$ This is due to the
fact that the $r-$analytic continuations of the solutions of the
polyharmonic equations in the annulus $A_{a,b}\subset \mathbb{R}^{d}$ live
on the set $\mathcal{A}_{a,b}$ but not on the set $\widetilde{A}_{a,b}$ !
The point is that in the case of the annulus $\mathcal{A}_{a,b}$ we cannot
identify the point $z\theta $ with $\left( z,\theta \right) ,$ in other
words, $\left( z,\theta \right) $ is not identified with $\left( -z,-\theta
\right) .$ 
\end{remark}

The paper is organized as follows: in Section \ref{SclassicalHardyspaces} we
recall background material about the Hardy space $H^{2}\left( \mathbb{B}%
_{R}\right) .$ In Section \ref{SpolyHardy} we introduce the polyharmonic
Hardy space $H^{2}\left( \mathcal{B}_{b}\right) $ on the ball $\mathcal{B}$
of the Klein-Dirac quadric. We prove that it is a Hilbert space, a maximum
principle, and infinite-differentiability of the functions in $H^{2}\left( 
\mathcal{B}_{b}\right) .$ In Section \ref{SSzegoKernel4HK} we construct a
Cauchy type kernel for $H^{2}\left( \mathcal{B}_{b}\right) .$ In Section \ref%
{SmainProperties} we prove other main properties of the space ... which
generalize similar properties of the one-dimensional Hardy spaces. In
Section \ref{Sbvp} we characterize the polyharmonic functions which are
extendible to the Hardy space $H^{2}\left( \mathcal{B}_{b}\right) .$ In
Section \ref{Scubature} we prove some of the main results of the paper,
about the error estimate of the polyharmonic interpolation, and about the
polyharmonic interpolatory cubature formulas, generalizing the polyharmonic
Gau\ss -Jacobi cubature formulas introduced in \cite{kounchevRenderArxiv}, 
\cite{kounchevRenderArkiv}.


\section{Classical Hardy spaces -- a reminder \label{SclassicalHardyspaces}}

Hardy spaces are a bridge between Harmonic and Complex Analysis. This is
based on the fact that the Taylor coefficients of a function $f\left(
x\right) $ on the real line $\mathbb{R}$ are at the same time the
coefficients of an orthogonal expansion of the analytic continuation $%
f\left( z\right) $ with respect to the basis $\left\{ z^{j}\right\} _{j\geq
0}$ which is orthogonal on the circle. Thus in a certain sense the setting
of the Hardy spaces represents a study of the properties of the real
functions having Taylor coefficients with $\sum \left\vert a_{j}\right\vert
^{2}<\infty $ by the methods of Complex Analysis. Since we are generalizing
before all the Hardy space $H^{2}$ we will recall the main results about it.
Let us put 
\begin{equation*}
M\left( f;r\right) :=\left\{ \frac{1}{2\pi }\int_{0}^{2\pi }\left\vert
f\left( re^{i\varphi }\right) \right\vert ^{2}d\varphi \right\} ^{1/2}.
\end{equation*}%
Then the for every function $f$ which is analytic in the disc $\mathbb{D}$
we define the Hardy space norm 
\begin{equation*}
\left\Vert f\right\Vert _{H^{2}\left( \mathbb{D}\right) }:=\sup_{r<1}M\left(
f;r\right) .
\end{equation*}

\begin{remark}
For every analytic function $f$ the function $M\left( f;r\right) $ is an
increasing function of $r,$ cf. Theorem $17.6$ in \cite{rudin}.
\end{remark}

A basic fact is that $H^{2}$ is a  Hilbert space and may be identified with
the limit values on the circle which coincide with $L^{2}\left( \mathbb{S}%
^{1}\right) .$ For every $g\in L^{2}\left( \mathbb{S}^{1}\right) $ we have
the norm defined by 
\begin{equation*}
\left\Vert g\right\Vert _{L_{2}\left( \mathbb{S}\right) }^{2}=\frac{1}{2\pi }%
\int_{-\pi }^{\pi }\left\vert g\left( e^{i\varphi }\right) \right\vert
^{2}d\varphi ,
\end{equation*}%
and the Fourier coefficients 
\begin{equation*}
\widehat{g}\left( n\right) =\frac{1}{2\pi }\int_{-\pi }^{\pi }g\left(
e^{i\varphi }\right) e^{-i\pi \varphi }d\varphi \qquad \text{for }n\in 
\mathbb{Z}.
\end{equation*}%
Some of the main properties of the classical Hardy spaces $H^{2}$ are
summarized in the following theorem (cf. Theorem $17.10$ in \cite{rudin}).

\begin{theorem}
\label{TclassicalHardy2} 1. An analytic function $f$ on $\mathbb{D}$ of the
type 
\begin{equation}
f\left( z\right) =\sum_{j=0}^{\infty }f_{j}z^{j}\qquad \text{for }z\in 
\mathbb{D}  \label{fzTaylor}
\end{equation}%
belongs to $H^{2}\left( \mathbb{D}\right) $ if and only if 
\begin{equation}
\sum_{j=0}^{\infty }\left\vert f_{j}\right\vert ^{2}<\infty ;
\label{sum_aj_lss_infty}
\end{equation}%
in that case 
\begin{equation*}
\left\Vert f\right\Vert _{H^{2}}^{2}=\sum_{j=0}^{\infty }\left\vert
f_{j}\right\vert ^{2}.
\end{equation*}

2. If $f\in H^{2}\left( \mathbb{D}\right) $ then $f$ has radial limits $%
f^{\ast}\left( e^{i\varphi}\right) $ at almost all points on the circle $%
\mathbb{S}$ and $f^{\ast}\in L^{2}\left( \mathbb{S}\right) .$ The Riesz
condition holds, i.e. 
\begin{align}
f_{j}^{\ast} & =0\qquad\text{for all }j<0;  \label{RieszConditionClassical}
\\
f_{j}^{\ast} & =f_{j}\qquad\text{for all }\geq0.  \notag
\end{align}
The $L^{2}-$approximation holds 
\begin{equation*}
\lim_{r\rightarrow1}\frac{1}{2\pi}\int_{0}^{2\pi}\left| f\left(
re^{i\varphi}\right) -f^{\ast}\left( e^{i\varphi}\right) \right|
^{2}d\varphi=0.
\end{equation*}
The integral of Poisson and of Cauchy of $f^{\ast}$ recover $f,$ i.e. 
\begin{align*}
f\left( z\right) & =\frac{1}{2\pi}\int_{0}^{2\pi}P_{r}\left(
\varphi-t\right) f^{\ast}\left( e^{it}\right) dt \\
f\left( z\right) & =\frac{1}{2\pi i}\int_{\Gamma}\frac{f^{\ast}\left(
\zeta\right) }{\zeta-z}d\zeta
\end{align*}
where $\Gamma$ is the positively oriented circle $\mathbb{S}.$

3. The mapping $f\longmapsto f^{\ast}$ is an isometry of $H^{2}$ on the
subspace of $L^{2}\left( \mathbb{S}\right) $ which consists of those $g$ for
which $\widehat{g}\left( j\right) =0$ for all $j<0.$
\end{theorem}

Let us recall the famous theorem of brothers F. and M. Riesz which concludes
the absolute continuity of a Borel measure on $\mathbb{S}^{1}$ only from the
annihilation of half of its Fourier coefficients, \cite[Theorem 17.13]{rudin}%
.

\begin{theorem}
\label{TFandMRiesz} Let $\mu $ be a complex valued Borel measure on the
circle $\mathbb{S}.$ If 
\begin{equation}
\int_{0}^{2\pi }e^{imt}d\mu \left( t\right) =0\qquad \text{for }m\geq 1,
\label{RieszConditions}
\end{equation}%
then the measure $\mu $ is absolutely continuous with respect to the
Lebesgue measure, i.e. there exists a function $f^{\ast }\in L^{1}\left( 
\mathbb{S}\right) $ such that $d\mu \left( t\right) =\frac{1}{2\pi }f^{\ast
}\left( t\right) dt$ for $t\in \left[ 0,2\pi \right] .$
\end{theorem}

\section{The polyharmonic Hardy space $H^{2}\left( \mathcal{B}\right) $ on
the ball of the Klein-Dirac quadric \label{SpolyHardy}}


At first we observe that the basis functions $b_{k,\ell ;j}\left( r,\theta
\right) $  defined in (\ref{PolyhBasis}) have a natural $r-$analytic
extension 
\begin{equation}
b_{k,\ell ,j}\left( z,\theta \right) :=z^{2j+k}Y_{k,\ell }\left( \theta
\right) .  \label{PolyhBasis1}
\end{equation}

Let us present some heuristics for explaining our main goal: We want to
define the Hardy space $H^{2}\left( \mathcal{B}_{R}\right) $ as a space of
functions which are uniform limits of sequences of complexified polynomials $%
P\left( z\theta \right) $ on compacts of the ball $\mathcal{B}_{R}$ in $%
\func{KDQ}.$ For that reason we need an appropriate inner product. Let us
make the \textbf{important observation} that there is a natural \emph{inner
product} where the basic functions (\ref{PolyhBasis1}) are \textbf{orthogonal%
}. Hence, for functions defined on the set $\mathbb{S}^{1}\times \mathbb{S}%
^{d-1}$ we introduce the inner product 
\begin{equation}
\left\langle f,g\right\rangle _{\ast }:=\frac{1}{2\pi }\int_{\mathbb{S}%
^{d-1}}\int_{0}^{2\pi }f\left( e^{i\varphi },\theta \right) \overline{%
g\left( e^{i\varphi },\theta \right) }d\varphi d\theta .
\label{scalarHardyBall}
\end{equation}%
The \textbf{crux} of our approach is the orthogonality of the basis
functions $b_{k,\ell ,j}$ in (\ref{PolyhBasis}) on the boundary of the
Klein-Dirac quadric $\mathbb{S}^{1}\times \mathbb{S}^{d-1}/\mathbb{Z}_{2},$
i.e. 
\begin{equation}
\left\langle b_{k,\ell ,j},b_{k^{\prime },\ell ,j}\right\rangle =\delta
_{k,k^{\prime }}\delta _{\ell ,\ell ^{\prime }}\delta _{j,j^{\prime }},
\label{bkljOrthogonal}
\end{equation}%
where the Kronecker symbol $\delta $ means $\delta _{\alpha ,\beta }=1$ for $%
\alpha =\beta $ and $0$ for $\alpha \neq \beta .$ This property is a
remarkable generalization of the orthogonality of the basis $\left\{
z^{j}\right\} _{j\geq 0}$ on the circle $\mathbb{S}^{1}$ and traces the
analogy to the one-dimensional case.

Further, we provide some arguments about the proper definition of the norm
of the prospective Hardy space. The objects of our polyharmonic Hardy space $%
H^{2}\left( \mathcal{B}\right) $ will be functions $f\left( z,\theta \right) 
$ which are representable as infinite sums in the $L^{2}$ sense, and are
absolutely and uniformly convergent on compacts with $\left\vert
z\right\vert <1$: 
\begin{align}
f\left( z,\theta \right) & =\sum_{k=0}^{\infty }\sum_{\ell =1}^{a_{k}}\left(
\sum_{j=0}^{\infty }f_{k,\ell ,j}z^{2j}\right) z^{k}Y_{k,\ell }\left( \theta
\right) .  \label{fAlmansiNEW} \\
& =\sum_{k=0}^{\infty }\sum_{\ell =1}^{a_{k}}f_{k,\ell }\left( z\right)
Y_{k,\ell }\left( \theta \right)   \label{fzthetaSeries}
\end{align}%
If $f_{k,\ell ,j}=0$ for all $j\geq N$ and $k\leq K$ for some $K\geq 0,$
then the function $f$ is a polynomial satisfying $\Delta ^{N}f\left(
x\right) =0$ for all $x\in \mathbb{R}^{d}.$ Hence, formula (\ref{fAlmansiNEW}%
) represents the $r-$complexification for the polynomials. 

In view of representation (\ref{fAlmansiNEW}) we introduce the following
naturally defined subspaces of the classical Hardy spaces on the unit disc $%
\mathbb{D}\subset \mathbb{C}$: 

\begin{definition}
The ''component spaces'' $H^{2,k}\left( \mathbb{D}\right) \subset
H^{2}\left( \mathbb{D}\right) $ consist of functions $f\in H^{2}\left( 
\mathbb{D}\right) $ having the representation $f\left( z\right) =f_{1}\left(
z^{2}\right) z^{k}.$ We put 
\begin{equation}
H^{2,k}\left( \mathbb{D}\right) :=\left\{ f\left( z\right) :f\left( z\right)
=f_{1}\left( z^{2}\right) z^{k},\ f_{1}\in H^{2}\left( \mathbb{D}\right)
\right\} .  \label{H2k}
\end{equation}
\end{definition}

Thus the space $H^{2,k}\left( \mathbb{D}\right) $ consists of the analytic
functions in $H^{2}\left( \mathbb{D}\right) $ having Taylor series 
\begin{equation*}
f\left( z\right) =\sum_{j=0}^{\infty }a_{j}z^{k+2j},
\end{equation*}%
respectively the norm on $H^{2,k}\left( \mathbb{D}\right) $ is the inherited
from $H^{2}\left( \mathbb{D}\right) .$ 

\begin{definition}
\label{DHardyball} We define the polyharmonic Hardy space $H^{2}\left( 
\mathcal{B}\right) $ on the unit ball $\mathcal{B}=\mathcal{B}_{1}$ of the
Klein-Dirac quadric defined in (\ref{BallComplexified}), as the space of
functions $f$ given by the Laplace-Fourier series (\ref{fAlmansiNEW}) with
coefficients $f_{k,\ell }\in H^{2,k}\left( \mathbb{D}\right) $ satisfying 
\begin{equation*}
\left\Vert f\right\Vert :=\sqrt{\sum_{k=0}^{\infty }\sum_{\ell
=1}^{a_{k}}\left\Vert f_{k,\ell }\right\Vert _{H^{2}\left( \mathbb{D}\right)
}^{2}}<\infty .
\end{equation*}
\end{definition}

\begin{remark}
The reader may note that Definition \ref{DHardyball} mimics the definition
of the classical Hardy spaces which are obtained as the closure of the
polynomials.
\end{remark}

Now let $f\in H^{2}\left( \mathcal{B}\right) .$ By  Definition \ref%
{DHardyball},  since all $f_{k,\ell }\in H^{2}\left( \mathbb{D}\right) $ it
follows that for $r\rightarrow 1^{-}$ and $z=re^{i\varphi }$ all $f_{k,\ell
}\left( re^{i\varphi }\right) $ have limiting values $f_{k,\ell }^{\ast
}\left( e^{i\varphi }\right) $ in $L_{2}\left( \mathbb{S}^{1}\right) ,$
hence 
\begin{equation}
\sum_{k=0}^{\infty }\sum_{\ell =1}^{a_{k}}\left\Vert f_{k,\ell }\right\Vert
_{H^{2}\left( \mathbb{D}\right) }^{2}=\sum_{k=0}^{\infty }\sum_{\ell
=1}^{a_{k}}\int_{0}^{2\pi }\left\vert f_{k,\ell }^{\ast }\left( e^{i\varphi
}\right) \right\vert ^{2}d\varphi =\sum_{k=0}^{\infty }\sum_{\ell
=1}^{a_{k}}\int_{0}^{2\pi }\left\vert f_{k,\ell }^{\ast }\left( z^{2}\right)
z^{k}\right\vert ^{2}d\varphi <\infty .  \label{fnormHKnice}
\end{equation}%
This implies 
\begin{equation*}
\lim_{r\rightarrow 1^{-}}\dint_{0}^{2\pi }\dint_{\mathbb{S}^{d-1}}\left\vert
f\left( re^{i\varphi },\theta \right) \right\vert ^{2}d\varphi d\theta
=\sum_{k=0}^{\infty }\sum_{\ell =1}^{a_{k}}\dint_{0}^{2\pi }\left\vert
f_{k,\ell }^{\ast }\left( e^{i\varphi }\right) \right\vert ^{2}d\varphi
<\infty ,
\end{equation*}%
hence,  it follows that $f\in L_{2}\left( \mathbb{D}\times \mathbb{S}%
^{d-1}\right) ,$ and due to (\ref{fAlmansiNEW}), also $f\in L_{2}\left( 
\mathcal{B}\right) .$ 

Note that as in the classical Hardy spaces the inner product (\ref%
{scalarHardyBall}) is good only for the polynomials but might not be well
defined for arbitrary functions having bad boundary behavior. For that
reason we have to change the definition of this inner product.

\begin{definition}
We put 
\begin{equation}
\left\langle f,g\right\rangle _{H^{2}\left( \mathcal{B}\right) }:=\frac{1}{%
2\pi }\lim_{\substack{ r\rightarrow 1 \\ r<1}}\int_{\mathbb{S}%
^{d-1}}\int_{0}^{2\pi }f\left( re^{i\varphi },\theta \right) \overline{%
g\left( re^{i\varphi },\theta \right) }d\varphi d\theta .
\label{scalarHardyBallGOOD}
\end{equation}
\end{definition}

The following theorem justifies our arguments above and is an analog to
results for the classical Hardy spaces, cf. Theorem \ref{TclassicalHardy2}
(or Theorem $17.10$ in \cite{rudin}).

\begin{theorem}
\label{THKbasic} 1. The space $H^{2}\left( \mathcal{B}\right) $ is complete.

2. It coincides with the space of functions $f$ having representation (\ref%
{fAlmansiNEW}) 
\begin{equation}
f\left( z,\theta \right) =\sum_{k=0}^{\infty }\sum_{\ell
=1}^{a_{k}}\sum_{j=0}^{\infty }f_{k,\ell ;j}z^{k+2j}Y_{k,\ell }\left( \theta
\right)   \label{fztheta}
\end{equation}%
with coefficients satisfying 
\begin{equation}
\sum_{k=0}^{\infty }\sum_{\ell =1}^{a_{k}}\sum_{j=0}^{\infty }\left\vert
f_{k,\ell ;j}\right\vert ^{2}<\infty .  \label{fnormHKsequence0}
\end{equation}

3. The norm of $f$ is given by 
\begin{equation*}
\left\Vert f\right\Vert _{H^{2}\left( \mathcal{B}\right) }=\left\{
\sum_{k=0}^{\infty }\sum_{\ell =1}^{a_{k}}\sum_{j=0}^{\infty }\left\vert
f_{k,\ell ;j}\right\vert ^{2}\right\} ^{1/2}.
\end{equation*}%
The following equality holds, for $z=re^{i\varphi },$ 
\begin{equation}
\left\Vert f\right\Vert _{H^{2}\left( \mathcal{B}\right)
}^{2}=\sum_{k=0}^{\infty }\sum_{\ell =1}^{a_{k}}\lim_{r\rightarrow 1}\frac{1%
}{2\pi }\int_{0}^{2\pi }\left\vert z^{k}\sum_{j=0}^{\infty }f_{k,\ell
,j}z^{2j}\right\vert ^{2}d\varphi =\sum_{k=0}^{\infty }\sum_{\ell
=1}^{a_{k}}\sum_{j=0}^{\infty }\left\vert f_{k,\ell ,j}\right\vert ^{2},
\label{fnormHKsequence}
\end{equation}

4. Every element $f\in $ $H^{2}\left( \mathcal{B}\right) $ is the limit of a
sequence of polynomials $P_{N}\in \mathcal{P}$ $\ $which satisfy 
\begin{equation*}
\Delta ^{N}P_{N}\left( x\right) =0\qquad \text{for }x\in B\subset \mathbb{R}%
^{d}.
\end{equation*}
\end{theorem}

We will prove only some of the above statements.

\proof
1. It is obvious by a direct estimation that condition $f\in $ $H^{2}\left( 
\mathcal{B}\right) $ implies the absolute convergence of the series (\ref%
{fztheta}) on every compact $K\subset \mathbb{D\times S}^{d-1}.$ Hence,
every function defined by the series (\ref{fztheta}) and satisfying (\ref%
{fnormHKsequence0}) also satisfies 
\begin{equation*}
\sum_{k=0}^{\infty }\sum_{\ell =1}^{a_{k}}\left\Vert f_{k,\ell }\right\Vert
_{H^{2}\left( \mathbb{B}_{1}\right) }^{2}=\lim_{r\rightarrow
1}\sum_{k=0}^{\infty }\sum_{\ell =1}^{a_{k}}\sum_{j=0}^{\infty }\left\vert
f_{k,\ell ;j}\right\vert ^{2}r^{2\left( k+2j\right) }=\sum_{k=0}^{\infty
}\sum_{\ell =1}^{a_{k}}\sum_{j=0}^{\infty }\left\vert f_{k,\ell
;j}\right\vert ^{2}<\infty ,
\end{equation*}%
which shows that the space $H^{2}\left( \mathcal{B}\right) $ is a tensor sum
of the $\ell _{2}-$spaces $\left\{ f_{k,\ell ;j}\right\} _{j\geq 0}$ taken
over all indices $\left( k,\ell \right) .$ Now, in a standard manner, the
completeness of the space $H^{2}\left( \mathcal{B}\right) $ follows from the
completeness of the space of sequences $\ell ^{2}.$

4. We have to see that every function defined by the series (\ref{fztheta})
and satisfying (\ref{fnormHKsequence0}) is a limit in the norm of $%
H^{2}\left( \mathcal{B}\right) $ of a sequence of polynomials. Indeed, for
every $\varepsilon >0$ we may find integers $k_{1}>0$ and $N>0$ such that 
\begin{equation*}
\sum_{k=k_{1}+1}^{\infty }\sum_{\ell =1}^{d_{k}}\sum_{j=N+1}^{\infty
}\left\vert f_{k,\ell ;j}\right\vert ^{2}<\varepsilon .
\end{equation*}%
Then the polynomial $P$ defined by 
\begin{equation*}
P\left( x\right) :=\sum_{k=0}^{k_{1}}\sum_{\ell
=1}^{d_{k}}\sum_{j=0}^{N}f_{k,\ell ;j}r^{k+2j}Y_{k,\ell }\left( \theta
\right) 
\end{equation*}%
satisfies with $z=re^{i\varphi }$ the following: 
\begin{eqnarray*}
\frac{1}{2\pi }\dint_{\mathbb{S}^{d-1}}\dint_{0}^{2\pi }\left\vert f\left(
z\theta \right) -P\left( z\theta \right) \right\vert ^{2}d\varphi d\theta 
&=&\sum_{k=k_{1}+1}^{\infty }\sum_{\ell =1}^{d_{k}}\sum_{j=N+1}^{\infty
}\left\vert f_{k,\ell ;j}\right\vert ^{2}r^{2\left( k+2j\right) } \\
&\leq &\sum_{k=k_{1}+1}^{\infty }\sum_{\ell =1}^{d_{k}}\sum_{j=N+1}^{\infty
}\left\vert f_{k,\ell ;j}\right\vert ^{2}\leq \varepsilon .
\end{eqnarray*}%
This implies 
\begin{equation*}
\left\Vert f-P\right\Vert _{H^{2}\left( \mathcal{B}\right) }\leq \varepsilon
,
\end{equation*}%
which ends the proof. 
\endproof%

\begin{remark}
The essence of Theorem \ref{THKbasic} is that, as in the classical Hardy
spaces, only the information about $f$ in the real domain, provided by the
Laplace-Fourier coefficients 
\begin{equation*}
f\left( x\right) =\sum_{k=0}^{\infty }\sum_{\ell =1}^{a_{k}}f_{k,\ell
}\left( r^{2}\right) r^{k}Y_{k,\ell }\left( \theta \right)
=\sum_{k=0}^{\infty }\sum_{\ell =1}^{a_{k}}\sum_{j=0}^{\infty }f_{k,\ell
;j}r^{k+2j}Y_{k,\ell }\left( \theta \right) ,
\end{equation*}%
determine when does $f$ belong to $H^{2}\left( \mathcal{B}\right) .$ All we
need to know is that they satisfy the convergence condition (\ref%
{fnormHKsequence0}), $\sum_{k=0}^{\infty }\sum_{\ell
=1}^{a_{k}}\sum_{j=0}^{\infty }\left\vert f_{k,\ell ,j}\right\vert
^{2}<\infty .$
\end{remark}

Theorem \ref{THKbasic} suggests the representation 
\begin{equation*}
H^{2}\left( \mathcal{B}\right) =\bigoplus_{k,\ell }^{\prime }H^{2,k},
\end{equation*}%
where the tensor sum is understood in the sense of equality (\ref%
{fzthetaSeries}), and the prime in the symbol $\bigoplus_{k,\ell }^{\prime }$
means that only sums are taken which are convergent as  (\ref%
{fnormHKsequence0})-(\ref{fnormHKsequence}).

The following result follows from the representation in (\ref{H2k}).

\begin{proposition}
For every function $f\in H^{2,k}\left( \mathbb{D}\right) $ the following
formula of Cauchy type holds: 
\begin{equation*}
f\left( \zeta\right) =\frac{1}{2\pi i}\int_{\Gamma_{1}}\frac{z}{%
z^{2}-\zeta^{2}}\frac{\zeta^{k}}{z^{k}}f\left( z\right) dz;
\end{equation*}
here $\zeta\in\mathbb{D}.$
\end{proposition}

\subsection{Cauchy type kernel for $H^{2}\left( \mathcal{B}\right) $ and
Hua-Aronszajn type formula \label{SSzegoKernel4HK}}

The orthogonality of the basis $\left\{ b_{k,\ell ,j}\right\} $ hints us to
construct a \textbf{Cauchy type kernel }and a corresponding \textbf{formula}
which reproduces the multivariate polynomials by using their values on the
set $\mathbb{S}^{1}\times \mathbb{S}^{d-1}/\mathbb{Z}_{2}.$ By exploiting
the above orthogonality, we obtain such formula easily, following the
general principles of constructing kernels, by putting: 
\begin{align}
K\left( z^{\prime },\theta ^{\prime };z,\theta \right) & :=\sum_{j\geq
0}\sum_{k=0}^{\infty }\sum_{\ell =1}^{a_{k}}b_{k,\ell ,j}\left( \zeta
;\theta ^{\prime }\right) b_{k,\ell ,j}\left( z;\theta \right) 
\label{CauchyTypeKernel} \\
& =\sum_{j\geq 0}\sum_{k=0}^{\infty }\sum_{\ell =1}^{a_{k}}\zeta
^{2j+k}Y_{k,\ell }\left( \theta ^{\prime }\right) z^{2j+k}Y_{k,\ell }\left(
\theta \right) .  \notag
\end{align}%
Note that this kernel is \textbf{absolutely} \textbf{convergent} for every $%
\zeta $ and $\theta $ with $\left\vert \zeta \right\vert <1,$ $\theta
^{\prime }\in \mathbb{S}^{d-1},$ and $z=e^{i\varphi },$ due to the estimates
for the spherical harmonics (\ref{SphericalEstimate}). By the orthogonality
property (\ref{bkljOrthogonal}), for every polynomial $P$ we obtain the
Cauchy type formula formula 
\begin{align}
P\left( \zeta \theta ^{\prime }\right) & =\left\langle K\left( \zeta ,\theta
^{\prime };\cdot \right) ,P\left( \cdot \right) \right\rangle _{\ast }
\label{CauchyTypeFormula} \\
& =\frac{1}{2\pi }\int_{\mathbb{S}^{d-1}}\int_{0}^{2\pi }K\left( \zeta
,\theta ^{\prime };z,\theta \right) \overline{P\left( z,\theta \right) }%
d\varphi d\theta .  \notag
\end{align}%
Let us remark that no such formula is available in the real domain. Formula (%
\ref{CauchyTypeFormula}) is a strong motivation to consider further the
consequences of the inner product (\ref{scalarHardyBall}).

Obviously, 
\begin{align*}
K\left( \zeta ,\theta ^{\prime };z,\theta \right) & =\sum_{j\geq 0}\left(
\zeta z\right) ^{2j}\sum_{k=0}^{\infty }\sum_{\ell =1}^{a_{k}}\left( \zeta
z\right) ^{k}Y_{k,\ell }\left( \theta ^{\prime }\right) Y_{k,\ell }\left(
\theta \right)  \\
& =\frac{1}{1-\zeta ^{2}z^{2}}\sum_{k=0}^{\infty }\sum_{\ell
=1}^{a_{k}}\left( \zeta z\right) ^{k}Y_{k,\ell }\left( \theta ^{\prime
}\right) Y_{k,\ell }\left( \theta \right) 
\end{align*}%
which shows the absolute convergence for $\left\vert \zeta z\right\vert <1.$
Let us recall that the usual \textbf{Poisson kernel} (see \cite{steinWeiss},
chapter $2,$ Theorem $1.9$) is given by 
\begin{equation*}
K_{P}\left( r,\theta ,\theta ^{\prime }\right) :=\frac{1-r^{2}}{\left\vert
\theta -r\theta ^{\prime }\right\vert ^{d}}=\sum_{k=0}^{\infty }\sum_{\ell
=1}^{d_{k}}r^{k}Y_{k,\ell }\left( \theta \right) Y_{k,\ell }\left( \theta
^{\prime }\right) \qquad \text{for }r<1
\end{equation*}%
in every dimension $d\geq 2.$ This expression is obviously close to the
above expression for $K\left( \zeta ,\theta ^{\prime };z,\theta \right) .$
We will apply the idea for the $r-$complexification to the kernel $%
K_{P}\left( r,\theta ,\theta ^{\prime }\right) $ and we will relate it to
the Cauchy type kernel $K\left( \zeta ,\theta ^{\prime };z,\theta \right) .$

\begin{proposition}
\label{PPoissonComplexified} For every complex number $w$ with $\left\vert
w\right\vert <1$ the following equality holds: 
\begin{equation}
\sum_{k=0}^{\infty }\sum_{\ell =1}^{d_{k}}w^{k}Y_{k,\ell }\left( \theta
\right) Y_{k,\ell }\left( \theta ^{\prime }\right) =\sum_{k=0}^{\infty
}w^{k}Z_{\theta }^{\left( k\right) }\left( \theta ^{\prime }\right) =\frac{%
1-w^{2}}{\left( 1-2w\left\langle \theta ,\theta ^{\prime }\right\rangle
+w^{2}\right) ^{\frac{d}{2}}},  \label{seriesPoissonComplexified}
\end{equation}%
where $Z_{\theta }^{\left( k\right) }\left( \theta ^{\prime }\right) $ are
the zonal harmonics, see \cite{steinWeiss}. Let us put $\cos \phi
=\left\langle \theta ,\theta ^{\prime }\right\rangle .$ For $d=2$ we have 
\begin{equation}
K_{P}\left( w,\theta ,\theta ^{\prime }\right) =\frac{1}{2\pi }+\frac{1}{\pi 
}\dsum_{k=1}^{\infty }w^{k}\frac{\cos k\phi }{\pi }=\frac{1}{2\pi }\frac{%
1-w^{2}}{1-2w\cos \phi +w^{2}}  \label{KP2d}
\end{equation}%
and for $d>2,$ 
\begin{equation}
K_{P}\left( w,\theta ,\theta ^{\prime }\right) =\dsum_{k=0}^{\infty
}w^{k}c_{k,d}P_{k}^{\lambda }\left( \cos \phi \right)   \label{KPnd}
\end{equation}%
where 
\begin{equation*}
c_{k,d}=\frac{1}{\omega _{d-1}}\frac{2k+d-2}{d-2},\qquad \lambda =\frac{d-2}{%
2}
\end{equation*}%
and $P_{k}^{\lambda }$ are the Legendre polynomials.

Hence, the $r-$complexification of the \textbf{Poisson kernel} $K_{P}\left(
r,\theta ,\theta ^{\prime }\right) $ is given by 
\begin{align}
K_{P}\left( w,\theta ,\theta ^{\prime }\right) & =\frac{1-w^{2}}{\left(
1-2w\left\langle \theta ,\theta ^{\prime }\right\rangle +w^{2}\right) ^{%
\frac{d}{2}}}  \label{PoissonKernelComplexified} \\
& =\frac{1-w^{2}}{\left( \left( 1-e^{i\phi }w\right) \left( 1-e^{-i\phi
}w\right) \right) ^{d/2}}.  \notag
\end{align}%
The \textbf{Cauchy type kernel} $K\left( \zeta ,\theta ^{\prime };z,\theta
\right) $ defined in (\ref{CauchyTypeKernel}) is related to the $r-$%
complexification of the Poisson kernel by the equality 
\begin{align}
K\left( \zeta ,\theta ^{\prime };z,\theta \right) & =\frac{1}{1-\zeta
^{2}z^{2}}K_{P}\left( \zeta z,\theta ,\theta ^{\prime }\right) 
\label{Szegoe-via-Poisson} \\
& =\frac{1}{\left( 1-2\zeta z\left\langle \theta ,\theta ^{\prime
}\right\rangle +\zeta ^{2}z^{2}\right) ^{\frac{d}{2}}}.  \notag
\end{align}
\end{proposition}

\proof
The proof uses the following results from \cite{steinWeiss}.

1. Formula (\ref{KP2d}) follows from formula (2.7) in chapter $4$ in \cite%
{steinWeiss}.

2. Formula (\ref{KPnd}) follows from Lemma $2.8,$ Theorem $2.10$ and Theorem 
$2.14$ in chapter $4$ in \cite{steinWeiss}.

By the estimates for the spherical harmonics in (\ref{SphericalEstimate}) it
follows that the left-hand side of (\ref{seriesPoissonComplexified}) is
absolutely and uniformly convergent on every compact $K\subset \mathbb{%
D\times S}^{d-1}\times \mathbb{S}^{d-1}.$ Representation (\ref%
{PoissonKernelComplexified}) follows from formula 
\begin{equation}
1-2w\cos \phi +w^{2}=\left( 1-e^{i\phi }w\right) \left( 1-e^{-i\phi
}w\right) .  \label{factorization}
\end{equation}%
\endproof%

The definition of the kernel $K\left( \zeta ,\theta ^{\prime };z,\theta
\right) $ shows that for all polynomials $P$ with real coefficients holds 
\begin{equation*}
P\left( \zeta \theta ^{\prime }\right) =\frac{1}{2\pi \omega _{d}}\int_{%
\mathbb{S}^{d-1}}\int_{0}^{2\pi }K\left( \zeta ,\theta ^{\prime };z,\theta
\right) \overline{P\left( z,\theta \right) }d\varphi d\theta \qquad \text{%
for }z=e^{i\varphi }.
\end{equation*}%
Since $\overline{z}=z^{-1},$ a change in the integration $\varphi
\rightarrow -\varphi $ shows that 
\begin{align*}
P\left( \zeta \theta ^{\prime }\right) & =\frac{1}{2\pi \omega _{d}}\int_{%
\mathbb{S}^{d-1}}\int_{0}^{2\pi }K\left( \zeta ,\theta ^{\prime };z,\theta
\right) \overline{P\left( z,\theta \right) }d\varphi d\theta  \\
& =\frac{1}{2\pi \omega _{d}}\int_{\mathbb{S}^{d-1}}\int_{0}^{2\pi }K\left(
\zeta ,\theta ^{\prime };\frac{1}{z},\theta \right) P\left( z,\theta \right)
d\varphi d\theta  \\
& =\frac{1}{2\pi i\omega _{d}}\int_{\mathbb{S}^{d-1}}\int_{\Gamma _{1}}\frac{%
1}{z}K\left( \zeta ,\theta ^{\prime };\frac{1}{z},\theta \right) P\left(
z,\theta \right) dzd\theta 
\end{align*}%
where $\Gamma _{1}$ is the positively oriented circle $\mathbb{S}^{1}.$ On
the other hand we see that 
\begin{align*}
\frac{1}{z}K\left( \zeta ,\theta ^{\prime };\frac{1}{z},\theta \right) & =%
\frac{1}{z}\frac{1}{\left( 1-2\frac{\zeta }{z}\left\langle \theta ,\theta
^{\prime }\right\rangle +\frac{\zeta ^{2}}{z^{2}}\right) ^{\frac{d}{2}}} \\
& =\frac{z^{d-1}}{\left( z^{2}-2\zeta z\left\langle \theta ,\theta ^{\prime
}\right\rangle +\zeta ^{2}\right) ^{\frac{d}{2}}}.
\end{align*}%
The last by definition is up to a factor the Hua-Aronszajn kernel, see \cite%
{aron}, p. $126.$

The above identities motivate the definition of the well-known \textbf{%
Hua-Aronszajn kernel} $H\left( \zeta ,\theta ^{\prime };z,\theta \right) $ 
given by 
\begin{equation}
H\left( \zeta ,\theta ^{\prime };z,\theta \right) =\frac{1}{\omega _{d}}%
\frac{z^{d-1}}{\left( \zeta ^{2}-2\zeta z\left\langle \theta ,\theta
^{\prime }\right\rangle +z^{2}\right) ^{\frac{d}{2}}},
\label{Hua-AronszajnKernelBALL}
\end{equation}%
where $\omega _{d}=\pi ^{d/2}/\Gamma \left( d/2\right) $ is the surface of
the sphere (see \cite{aron}, p. $122$ Theorem $1.1,$ and p. $126,$ Remark $%
1.4,$ where up to a factor it is called \emph{Cauchy kernel} for the \emph{%
Cartan classical} domain $\mathcal{R}_{IV}$).

From above we see that the following equality holds:  
\begin{align}
H\left( \zeta ,\theta ^{\prime };z,\theta \right) & =\frac{1}{\omega _{d}}%
\frac{1}{z}K\left( \zeta ,\theta ^{\prime };\frac{1}{z},\theta \right) =%
\frac{1}{z}\frac{1}{1-\zeta ^{2}/z^{2}}\sum_{k=0}^{\infty }\sum_{\ell
=1}^{a_{k}}\left( \zeta /z\right) ^{k}Y_{k,\ell }\left( \theta ^{\prime
}\right) Y_{k,\ell }\left( \theta \right) 
\label{Hua-AronszajnKernelBALLseries} \\
& =\frac{1}{\omega _{d}}\frac{1}{z}\frac{1}{\left( 1-2\zeta /z\left\langle
\theta ,\theta ^{\prime }\right\rangle +\zeta ^{2}/z^{2}\right) ^{\frac{d}{2}%
}}.  \notag
\end{align}

\section{Main properties of the Hardy spaces $H^{2}\left( \mathcal{B}\right) 
$ \label{SmainProperties}}

In the next theorem we provide a generalization of the classical boundary
value properties of the Hardy spaces $H^{2},$ see e.g. Theorem \ref%
{TclassicalHardy2} (or Theorem $17.10,$ $17.12$ and $17.13$ in \cite{rudin}%
).  

\begin{theorem}
\label{THKbasic2} Let $f\in H^{2}\left( \mathcal{B}\right) .$

1. \textbf{Fatou type theorem}: For $r\rightarrow 1^{-},$ and for almost all 
$\varphi \in \left[ 0,2\pi \right] $ and $\theta \in \mathbb{S}^{d-1},$ the
function $f\left( re^{i\varphi }\theta \right) $ has a radial limit which we
denote by $f^{\ast }\left( e^{i\varphi },\theta \right) ,$ and which
satisfies $f^{\ast }\left( e^{i\varphi },\theta \right) \in L^{2}\left( 
\mathbb{S}\times \mathbb{S}^{d-1}\right) .$ If the expansion of the function 
$f$ is given by (\ref{fztheta}), 
\begin{equation*}
f\left( z,\theta \right) =\sum_{k=0}^{\infty }\sum_{\ell
=1}^{a_{k}}\sum_{j=0}^{\infty }f_{k,\ell ;j}z^{k+2j}Y_{k,\ell }\left( \theta
\right) \qquad \text{for all }\left\vert z\right\vert <1,\theta \in \mathbb{S%
}^{d-1},
\end{equation*}%
then $f^{\ast }\left( e^{i\varphi },\theta \right) $ is given by the
Laplace-Fourier series 
\begin{equation}
f^{\ast }\left( e^{i\varphi },\theta \right) =\sum_{k=0}^{\infty }\sum_{\ell
=1}^{a_{k}}\sum_{j=0}^{\infty }f_{k,\ell ;j}z^{k+2j}Y_{k,\ell }\left( \theta
\right) \qquad \text{for }z=e^{i\varphi }.  \label{fstarLaplaceFourier}
\end{equation}

2. The following \textbf{limiting relation} holds, 
\begin{equation*}
\lim_{r\rightarrow 1}\frac{1}{2\pi }\int_{0}^{2\pi }\int_{\mathbb{S}%
^{d-1}}\left\vert f\left( re^{i\varphi }\theta \right) -f^{\ast }\left(
e^{i\varphi }\theta \right) \right\vert ^{2}d\varphi d\theta =0.
\end{equation*}

3. Let the Laplace-Fourier series of the function $f^{\ast }\left(
e^{i\varphi },\theta \right) \in L^{2}\left( \mathbb{S}\times \mathbb{S}%
^{d-1}\right) $ be given by 
\begin{equation*}
f^{\ast }\left( e^{i\varphi },\theta \right) =\sum_{k=0}^{\infty }\sum_{\ell
=1}^{a_{k}}f_{k,\ell }^{\ast }\left( e^{i\varphi }\right) Y_{k,\ell }\left(
\theta \right) =\sum_{k=0}^{\infty }\sum_{\ell =1}^{a_{k}}\sum_{j=0}^{\infty
}f_{k,\ell ;j}^{\ast }z^{k+2j}Y_{k,\ell }\left( \theta \right) \qquad \text{%
for }z=e^{i\varphi }.
\end{equation*}%
Then the Fourier coefficients of the functions $f_{k,\ell }^{\ast }\left(
e^{i\varphi }\right) $ satisfy the following zero conditions, which we call 
\textbf{Riesz type conditions}: 
\begin{equation}
f_{k,\ell ;j}^{\ast }=0\qquad \text{for all }j\neq k,k+2,k+4,...\ .
\label{RieszType}
\end{equation}%
These conditions are equivalent to (the usual form of Riesz conditions) 
\begin{equation}
\int_{0}^{2\pi }\int_{\mathbb{S}^{d-1}}f^{\ast }\left( e^{i\varphi },\theta
\right) e^{-ij\varphi }Y_{k,\ell }\left( \theta \right) d\varphi d\theta
=0\qquad \text{for all }j\neq k,k+2,k+4,...\ .  \label{RieszType2}
\end{equation}

4. The following \textbf{Cauchy-Hua-Aronszajn} formula holds 
\begin{align}
f\left( \zeta ,\theta ^{\prime }\right) & =\frac{1}{2\pi i}\int_{\Gamma
_{1}}\int_{\mathbb{S}^{d-1}}H\left( \zeta ,\theta ^{\prime };z,\theta
\right) f^{\ast }\left( z,\theta \right) dzd\theta 
\label{HuaAronszajnFormula} \\
& =\frac{1}{2\pi \omega _{d}}\int_{0}^{2\pi }\int_{\mathbb{S}^{d-1}}K\left(
\zeta ,\theta ^{\prime };\frac{1}{z},\theta \right) f^{\ast }\left( z,\theta
\right) d\varphi d\theta   \notag \\
& =\frac{1}{2\pi \omega _{d}}\int_{0}^{2\pi }P_{r^{2}}\left( 2\varphi
-2\varphi ^{\prime }\right) \int_{\mathbb{S}^{d-1}}K_{P}\left( \frac{\zeta }{%
z},\theta ,\theta ^{\prime }\right) f^{\ast }\left( z,\theta \right) d\theta
d\varphi ,  \notag
\end{align}%
where we use the notations $z=e^{i\varphi }$ and $\zeta =re^{i\varphi
^{\prime }},$ the kernel $K$ is the \textbf{Cauchy type kernel} (\ref%
{CauchyTypeKernel}), and $P_{r}\left( \varphi \right) $ is the usual
two-dimensional Poisson kernel. We call the kernel 
\begin{equation}
P_{r^{2}}\left( 2\varphi -2\varphi ^{\prime }\right) K_{P}\left( \frac{\zeta 
}{z},\theta ,\theta ^{\prime }\right)   \label{PoissonTypeKernel}
\end{equation}%
the \textbf{modified} \textbf{Poisson type kernel}. Here the contour $\Gamma
_{1}$ is the positively oriented circle $\mathbb{S}^{1},$ or a bigger
contour which encircles it.

5. \textbf{Dirichlet problem} with $L^{2}$ data: If a function $f^{\ast
}\left( e^{i\varphi },\theta \right) \in L^{2}\left( \mathbb{S}\times 
\mathbb{S}^{d-1}\right) $ satisfies the Riesz type conditions (\ref%
{RieszType}) then there exists an unique function $f\in H^{2}\left( \mathcal{%
B}\right) $ which has as a ''non-tangential limit''\ the function $f^{\ast
}\left( e^{i\varphi },\theta \right) $ in the sense 
\begin{equation*}
\lim_{r\longrightarrow 1}\int_{\mathbb{S}^{d-1}}\int_{0}^{2\pi }\left|
f\left( re^{i\varphi },\theta \right) -f^{\ast }\left( e^{i\varphi },\theta
\right) \right| ^{2}d\varphi d\theta =0.
\end{equation*}

6. Every function $f^{\ast }\in L^{2}\left( \mathbb{S}^{1}\times \mathbb{S}%
^{d-1}\right) $ which satisfies the Riesz tyep conditions (\ref{RieszType})
belongs to the space $L^{2}\left( \mathbb{S}^{1}\times \mathbb{S}^{d-1}/%
\mathbb{Z}_{2}\right) .$ The map $f\longrightarrow f^{\ast }$ is an \textbf{%
isometry} between $H^{2}\left( \mathcal{B}\right) $ and the subspace of $%
L^{2}\left( \mathbb{S}^{1}\times \mathbb{S}^{d-1}/\mathbb{Z}_{2}\right) .$

7. The norm is given by 
\begin{equation*}
\left\Vert f\right\Vert _{H^{2}\left( \mathcal{B}\right) }^{2}=\frac{1}{2\pi 
}\int_{\mathbb{S}^{d-1}}\int_{0}^{2\pi }\left\vert f^{\ast }\left(
e^{i\varphi }\theta ^{\prime }\right) \right\vert ^{2}d\varphi d\theta
^{\prime }.
\end{equation*}
\end{theorem}

\proof%
1. We will omit  the proof of this non-trivial but standardly proved
statement.

2. By Theorem \ref{THKbasic} the function $f$ is represented by a series (%
\ref{fztheta}) satisfying $\sum_{k=0}^{\infty }\sum_{\ell
=1}^{a_{k}}\sum_{j=0}^{\infty }\left\vert f_{k,\ell ;j}\right\vert
^{2}<\infty .$ Let us consider the function 
\begin{equation*}
f^{\ast }\left( e^{i\varphi },\theta \right) :=\sum_{k=0}^{\infty
}\sum_{\ell =1}^{a_{k}}\sum_{j=0}^{\infty }f_{k,\ell ;j}z^{k+2j}Y_{k,\ell
}\left( \theta \right) \qquad \text{for }z=e^{i\varphi }.
\end{equation*}%
Since the functions $\left\{ e^{ij\varphi }Y_{k,\ell }\left( \theta \right)
:j\in \mathbb{Z},\ \text{all }\left( k,\ell \right) \right\} $ form a basis
of the space $L^{2}\left( \mathbb{S}^{1}\times \mathbb{S}^{d-1}\right) ,$ it
follows that $f^{\ast }\in L^{2}\left( \mathbb{S}^{1}\times \mathbb{S}%
^{d-1}\right) .$ By Parseval's theorem we see that 
\begin{equation*}
\frac{1}{2\pi }\int_{0}^{2\pi }\int_{\mathbb{S}^{d-1}}\left\vert f\left(
\rho e^{i\varphi }\theta \right) -f^{\ast }\left( e^{i\varphi },\theta
\right) \right\vert ^{2}d\varphi d\theta =\sum_{k=0}^{\infty }\sum_{\ell
=1}^{a_{k}}\sum_{j=0}^{\infty }\left\vert f_{k,\ell ;j}\right\vert
^{2}\left( 1-\rho ^{k+2j}\right) ^{2}.
\end{equation*}%
It is easy to see that for $\rho \rightarrow 1^{-}$ the last tends to $0.$

The proof of 3.) is evident due to the expansion (\ref{fstarLaplaceFourier}).

4. From formula (\ref{Hua-AronszajnKernelBALLseries}) it follows that for
every function $g\left( e^{i\varphi },\theta \right) \in L^{2}\left( \mathbb{%
S}^{1}\times \mathbb{S}^{d-1}\right) $ holds 
\begin{gather*}
\frac{1}{2\pi i\omega _{d}}\int_{\Gamma _{1}}\int_{\mathbb{S}^{d-1}}\frac{%
z^{d-1}f^{\ast }\left( z,\theta \right) }{r\left( \zeta \theta ^{\prime
}-z\theta \right) ^{d}}dzd\theta =\qquad \qquad \qquad \qquad \qquad \qquad
\qquad \qquad \qquad  \\
=\frac{1}{2\pi \omega _{d}}\int_{0}^{2\pi }\int_{\mathbb{S}^{d-1}}K\left(
\zeta ,\theta ^{\prime };z,\theta \right) f^{\ast }\left( \overline{z}%
,\theta \right) d\varphi d\theta \qquad \text{for }z=e^{i\varphi }.
\end{gather*}%
From the series representation of $f$ and $f^{\ast },$ by the definition of
the kernel $K\left( \zeta ,\theta ^{\prime };z,\theta \right) $ it follows 
\begin{equation*}
f\left( \zeta \theta ^{\prime }\right) =\frac{1}{2\pi }\int_{0}^{2\pi }\int_{%
\mathbb{S}^{d-1}}K\left( \zeta ,\theta ^{\prime };z,\theta \right) f^{\ast
}\left( \overline{z},\theta \right) d\varphi d\theta .
\end{equation*}%
Further, since 
\begin{equation*}
P_{r}\left( \varphi -\varphi ^{\prime }\right) =\sum_{j=-\infty }^{\infty
}r^{\left\vert j\right\vert }e^{ij\left( \varphi ^{\prime }-\varphi \right) }
\end{equation*}%
and 
\begin{equation*}
K_{P}\left( \frac{\zeta }{z},\theta ,\theta ^{\prime }\right)
=\sum_{k=0}^{\infty }\sum_{\ell =1}^{a_{k}}r^{k}e^{ik\left( \varphi ^{\prime
}-\varphi \right) }Y_{k,\ell }\left( \theta \right) Y_{k,\ell }\left( \theta
^{\prime }\right) 
\end{equation*}%
the product $P_{r^{2}}\left( 2\varphi -2\varphi ^{\prime }\right)
K_{P}\left( \frac{\zeta }{z},\theta ,\theta ^{\prime }\right) $ has terms
which in the integral will vanish due to the Riesz conditions satisfied by $%
f^{\ast }$, and there will remain only the terms corresponding to $K\left(
\zeta ,\theta ^{\prime };z,\theta \right) ,$ which proves the formula.

5. Every function $f^{\ast }\left( e^{i\varphi },\theta \right) $ which
satisfies the Riesz type conditions (\ref{RieszType}) has a Laplace-Fourier
series of the type (\ref{fstarLaplaceFourier}), say 
\begin{equation*}
f^{\ast }\left( e^{i\varphi },\theta \right) =\sum_{k=0}^{\infty }\sum_{\ell
=1}^{a_{k}}\sum_{j=0}^{\infty }a_{k,\ell ;j}z^{k+2j}Y_{k,\ell }\left( \theta
\right) \qquad \text{for }z=e^{i\varphi }.
\end{equation*}%
Since $f^{\ast }\left( e^{i\varphi },\theta \right) \in L^{2}\left( \mathbb{S%
}^{1}\times \mathbb{S}^{d-1}\right) ,$ the coefficients satisfy $%
\sum_{k=0}^{\infty }\sum_{\ell =1}^{a_{k}}\sum_{j=0}^{\infty }\left\vert
a_{k,\ell ;j}\right\vert ^{2}<\infty ,$ hence, we may define the function 
\begin{equation*}
f\left( z\theta \right) :=\sum_{k=0}^{\infty }\sum_{\ell
=1}^{a_{k}}\sum_{j=0}^{\infty }a_{k,\ell ;j}z^{k+2j}Y_{k,\ell }\left( \theta
\right) ;
\end{equation*}%
the limiting property follows as in item $2).$

Finally, item 6.) follows from the above arguments.

\endproof%

Theorem \ref{THKbasic2} shows that we may solve Boundary Value Problems in
the spaces $H^{2}\left( \mathcal{B}\right) $ which is an essential advantage
over the situation with the holomorphic functions in $\mathbb{C}^{d}$ and
alternative definitions of Hardy space in several dimensions, see \cite%
{stein}, \cite{rudin70}, \cite{steinWeiss}, \cite{rudin80}.

By Theorem \ref{THKbasic2}  the following space: 
\begin{equation*}
\left\{ f^{\ast }\in L^{2}\left( \mathbb{S}\times \mathbb{S}^{d-1}\right) :%
\text{ }f^{\ast }\text{ satisfies the Riesz type condition (\ref{RieszType})}%
\right\} 
\end{equation*}%
is isomorphic to the space $H^{2}\left( \mathcal{B}\right) .$ The
Cauchy-Hua-Aronszajn formula (\ref{HuaAronszajnFormula}) generalizes the
Cauchy formula in $\mathbb{C}$ and the Poisson formula in $\mathbb{R}^{d}$
at the same time.

\begin{remark}
Let us formulate a \textbf{conjecture} about an analog to brothers' Riesz
theorem: Let the complex valued Borel measure $\mu \left( \varphi ,\theta
\right) $ be given on $\mathbb{S}^{1}\times \mathbb{S}^{d-1}$ with $\varphi
\in \left[ 0,2\pi \right] $ and $\theta \in \mathbb{S}^{d-1}.$ Assume that
for all indices $\left( k,\ell \right) $ holds 
\begin{equation}
\int_{0}^{2\pi }\int_{\mathbb{S}^{d-1}}\overline{z}^{j}Y_{k,\ell }\left(
\theta \right) d\mu \left( \varphi ,\theta \right) =0\qquad \text{for }j\neq
k,k+2,k+4,...  \label{BrothersRiesz}
\end{equation}%
Is the measure $\mu $ \textquotedblright absolutely
continuous\textquotedblright , i.e. does there exist a function $f^{\ast }$
which is in $L_{1}\left( \mathbb{S}^{1}\times \mathbb{S}^{d-1}\right) $ such
that $d\mu \left( \varphi ,\theta \right) =f^{\ast }\left( \varphi ,\theta
\right) d\varphi d\theta $ ? It seems that the answer in this form is
negative, but a positive answer needs some additional properties of the
measure $\mu .$ A thorough discussion to this question will be considered in 
\cite{kounchevRenderBook}. Let us remark that a genuine analog of the
brothers Riesz theorem is difficult to achieve for all approaches to Hardy
spaces, cf. \cite{bochner}, \cite{helsonLowdenslager}, \cite{rudin70}, \cite%
{stein}, \cite{steinWeiss}, \cite{rudin80}.
\end{remark}

\begin{remark}
Another \textbf{conjecture} about the boundary properties of the remarkable
Poisson type kernel in (\ref{HuaAronszajnFormula}) generalizing a classical
situation is the following:\ Assume that the function $g\left( e^{i\varphi
},\theta \right) $ belongs to $C\left( \mathbb{S}^{1}\times \mathbb{S}^{d-1}/%
\mathbb{Z}_{2}\right) $ and  satisfies the brothers Riesz conditions (\ref%
{RieszType2}). Then by means of formula (\ref{fstarLaplaceFourier}) (or
equivalently, by (\ref{HuaAronszajnFormula})) we may define a function $%
F_{g}\left( z,\theta \right) $ on the interior (for $\left\vert z\right\vert
=r<1$ ) of the ball $\mathcal{B}$ of the Klein-Dirac quadric. Let us put $%
F_{g}\left( e^{i\varphi },\theta \right) =g\left( e^{i\varphi },\theta
\right) $ for $r=1.$ We conjecture that the function $F_{g}$ is continuos on
the closure $\overline{\mathcal{B}}.$ Let us note that in the classical case
the Poisson kernel is used to prove similar statement, see chapter $2,$
Theorem $1.9$ in \cite{steinWeiss}. Here we expect that the modified Poisson
kernel (\ref{PoissonTypeKernel}) will of central importance for the
solution. 
\end{remark}

\subsection{Maximum principle}

In the classical case of the Hardy spaces, the maximum principle is
intimately related to the Cauchy formula in $\mathbb{C}$ or to the Poisson
formula in $\mathbb{R}^{d}$ (see the proof of the completeness of $H^{p}$ in
Remark $17.8$ in \cite{rudin}). A  weak form of maximum principle alows to
prove that the elements of $H^{2}$ are uniform limits of polynomials on
compact subsets of $\mathbb{D}\mathbf{.}$ Here we prove analog to this for
the polyharmonic Hardy space  $H^{2}\left( \mathcal{B}\right) .$ In the next
theorem we see that the explicit form for the Cauchy-Hua-Aronszajn kernel is
essential for proving a maximum principle.

\begin{theorem}
\label{TmaximumPrincipleHK} Let $f\in H^{2}\left( \mathcal{B}\right) .$ For
every $q$ with $0<q<1$ we have the following (weak) maximum principle holds: 
\begin{equation*}
\left\vert f\left( \zeta \theta \right) \right\vert \leq \left( 1-q\right)
^{-d}\left\Vert f\right\Vert _{H^{2}\left( \mathcal{B}\right) }\qquad \text{%
for all }\left\vert \zeta \right\vert \leq q,\ \theta \in \mathbb{S}^{d-1}.
\end{equation*}
More generally, for every mixed derivative $D_{\zeta ,\theta }^{\alpha }$
with respect to the variables $\zeta $ and $\theta $, we have the maximum
principle 
\begin{align*}
\left\vert D^{\alpha }f\left( \zeta \theta \right) \right\vert & \leq
C_{1}\times \left( 1-q\right) ^{-d-\left\vert \alpha \right\vert }\left[
\left( \frac{d}{2}\right) \left( \frac{d}{2}+1\right) \cdot \cdot \cdot
\left( \frac{d}{2}+\left\vert \alpha \right\vert \right) \right] \left\vert
\alpha \right\vert \left\Vert f\right\Vert _{H^{2}\left( \mathcal{B}\right) }
\\
\qquad \text{for all }\left\vert \zeta \right\vert & \leq q,\ \text{and }%
\theta \in \mathbb{S}^{d-1};
\end{align*}%
here the constant $C_{1}>0$ is independent of $\alpha .$ Respectively, for
real $\zeta =r$ with $x=r\theta $ this gives an estimate for $D_{x}^{\alpha
}f\left( x\right) .$
\end{theorem}

\proof%
Indeed, for $\left\vert \zeta \right\vert \leq q$ and $z=e^{i\varphi },$ if
we put $w=\zeta /z,$ by the Hua-Aronszajn formula (\ref{HuaAronszajnFormula}%
), it follows 
\begin{align*}
\left\vert f\left( \zeta \theta \right) \right\vert & \leq \frac{1}{2\pi }%
\int_{\mathbb{S}^{d-1}}\int_{0}^{2\pi }\left\vert K\left( \zeta ,\theta
^{\prime };\frac{1}{z},\theta \right) \right\vert \left\vert f^{\ast }\left(
z\theta ^{\prime }\right) \right\vert \left\vert dz\right\vert d\theta
^{\prime } \\
& \leq \left\{ \frac{1}{2\pi }\int_{\mathbb{S}^{d-1}}\int_{0}^{2\pi
}\left\vert \left( 1-2w\left\langle \theta ,\theta ^{\prime }\right\rangle
+w^{2}\right) ^{-\frac{d}{2}}\right\vert ^{2}d\varphi d\theta ^{\prime
}\right\} ^{1/2}\times  \\
& \times \left\{ \frac{1}{2\pi }\int_{\mathbb{S}^{d-1}}\int_{0}^{2\pi
}\left\vert f^{\ast }\left( e^{i\varphi }\theta ^{\prime }\right)
\right\vert ^{2}d\varphi d\theta ^{\prime }\right\} ^{1/2}.
\end{align*}%
Here we apply (\ref{factorization}) with $\cos \psi =\left\langle \theta
,\theta ^{\prime }\right\rangle ,$ which implies the inequality 
\begin{equation*}
\left\vert \left( 1-2w\left\langle \theta ,\theta ^{\prime }\right\rangle
+w^{2}\right) \right\vert =\left\vert \left( 1-e^{i\psi }w\right) \left(
1-e^{-i\psi }w\right) \right\vert \geq \left( 1-q\right) ^{2},
\end{equation*}%
hence, 
\begin{align*}
\left\vert \left( 1-2w\left\langle \theta ,\theta ^{\prime }\right\rangle
+w^{2}\right) ^{\frac{d}{2}}\right\vert & =\left\vert \left( 1-e^{i\psi
}w\right) \left( 1-e^{-i\psi }w\right) \right\vert ^{d} \\
& \geq \left( 1-q\right) ^{d}.
\end{align*}%
By Theorem \ref{THKbasic2} this implies 
\begin{equation*}
\left\vert f\left( \zeta \theta \right) \right\vert \leq \left( 1-q\right)
^{-d}\left\Vert f\right\Vert _{H^{2}\left( \mathcal{B}\right) },
\end{equation*}%
which ends the proof of the first part of our statement.

In the same way, by differentiating under the sign of the integral in the
Hua-Aronszajn formula, we  obtain the estimate for the derivatives, 
\begin{equation*}
D_{\zeta ,\theta }^{\alpha }\left( \left( 1-2\frac{\zeta }{z}\left\langle
\theta ,\theta ^{\prime }\right\rangle +\frac{\zeta ^{2}}{z^{2}}\right) ^{-%
\frac{d}{2}}\right) .
\end{equation*}%
By the Leibniz differentiation formula, the differentiation with respect to $%
\zeta $ and $\theta $ gives the following estimate 
\begin{gather*}
\left\vert D_{\zeta ,\theta }^{\alpha }\left( \left( 1-2\frac{\zeta }{z}%
\left\langle \theta ,\theta ^{\prime }\right\rangle +\frac{\zeta ^{2}}{z^{2}}%
\right) ^{-\frac{d}{2}}\right) \right\vert \leq \qquad \qquad \qquad \qquad
\qquad \qquad  \\
\qquad \qquad \qquad \qquad \leq C\left( \frac{1}{1-q}\right) ^{d+\left\vert
\alpha \right\vert }\left\vert \left( -\frac{d}{2}\right) \left( -\frac{d}{2}%
-1\right) \cdot \cdot \cdot \left( -\frac{d}{2}-\left\vert \alpha
\right\vert \right) \right\vert \times \left\vert \alpha \right\vert ,
\end{gather*}%
and this ends the proof.

\endproof%

We have the following immediate corollary about the regularity of the
functions in the space $H^{2}\left( \mathcal{B}\right) .$

\begin{corollary}
The functions in $H^{2}\left( \mathcal{B}\right) $ belong to $C^{\infty
}\left( \mathcal{B}\right) .$
\end{corollary}

The proof follows from the maximum principle in Theorem \ref%
{TmaximumPrincipleHK} since every $f\in H^{2}\left( \mathcal{B}\right) $ and
the derivatives of $f$ are uniform limits of a sequence of polynomials $%
P_{N}\left( z\theta \right) $ and the respecitve derivatives of $P_{N}\left(
z\theta \right) $ on every compact sets $K\times \mathbb{S}^{d-1}$ where the
compact $K\subset \mathbb{D}.$

\subsection{Real analytic functions and the space $H^{2}\left( \mathcal{B}%
\right) $}

As is the case with the classical analytic functions the elements of $%
H^{2}\left( \mathcal{B}\right) $ are real analytic which we prove in the
next theorem.

\begin{theorem}
\label{THKrealAnalytic} Let $F\in H^{2}\left( \mathcal{B}\right) .$ Then the
function $f\left( x\right) =F\left( r\theta \right) $ is real-analytic in
the ball $B\subset \mathbb{R}^{d}.$
\end{theorem}

\proof
We use the maximum principle: By a change of the variables, and a routine
estimation with the Leibnitz rule, we obtain 
\begin{equation*}
\left| D^{\alpha }f\left( x\right) \right| \leq C\max_{\left| \beta \right|
=\left| \alpha \right| }\left| D_{r,\theta }^{\beta }f\left( r\theta \right)
\right| .
\end{equation*}
The last is estimated by the maximum principle on compacts, by Theorem \ref%
{TmaximumPrincipleHK}, which gives 
\begin{equation*}
\left| D_{r,\theta }^{\beta }f\left( r\theta \right) \right| \leq C\frac{%
\left| \beta \right| !}{M_{q}^{\left| \beta \right| }}\leq C_{1}\frac{\beta !%
}{M_{q}^{\left| \beta \right| }}\qquad \text{for }r\leq q<1.
\end{equation*}
The statement of the theorem follows by an equivalent definition of
real-analytic function (cf. \cite{aron}, p. $17$).

\endproof%

\section{BVPs for the polyharmonic operator $\Delta ^{N}$ and the spaces $%
H^{2}\left( \mathcal{B}\right) $ \label{Sbvp}}

\subsection{The one-dimensional case}

Let us provide some heuristical observations from the one-dimensional case.

As we already said, the Hardy space setting provides interpretation to real
Taylor series $f\left( t\right) =\sum_{j=0}^{\infty }a_{j}t^{j}$ for which
we have $\sum_{j=0}^{\infty }\left\vert a_{j}\right\vert ^{2}<\infty $ as
the Fourier series $\sum_{j=0}^{\infty }a_{j}e^{ij\varphi }$ in the Hardy
space $H^{2}\left( \mathbb{D}\right) .$ Thus the Hardy spaces provide a
different way to encode the information of \textbf{real-domain data}. There
is an alternative way to encode the information, by means of  boundary data,
e.g. by means of 
\begin{equation*}
\left\{ \frac{d^{2j}f}{dt^{2j}}\left( -1\right) ,\frac{d^{2j}f}{dt^{2j}}%
\left( 1\right) \right\} _{j\geq 0}.
\end{equation*}%
It is more convenient to consider further differential operators defined on
the functions on $\mathbb{S}.$ 

Let $F$ be a function defined on $\mathbb{S}.$ Indeed, if we have a function 
$F$ which is representable as a series 
\begin{equation*}
F\left( \varphi \right) =\sum_{j=0}^{\infty }a_{j}e^{ij\varphi }
\end{equation*}%
then we obtain 
\begin{equation*}
DF\left( \varphi \right) =\sum_{j=1}^{\infty }ja_{j}e^{i\left( j-1\right)
\varphi },
\end{equation*}%
where we have put 
\begin{equation*}
DF\left( e^{i\varphi }\right) :=\frac{1}{ie^{i\varphi }}F^{\prime }\left(
e^{i\varphi }\right) .
\end{equation*}

We want to investigate for which sequences $\left\{ c_{j}\right\}
_{j=0}^{\infty }$ and $\left\{ d_{j}\right\} _{j=0}^{\infty }$ there exists
a function $F\in H^{2}\left( \mathbb{D}\right) $ such that 
\begin{equation*}
D^{2j}F\left( -1\right) =c_{j},\quad D^{2j}F\left( 1\right) =d_{j}\qquad 
\text{for }j=0,1,2,...\ .
\end{equation*}%
since $D^{2j}$ is the operator $\Delta ^{j}$ in the one-dimensional case.
The matrix of the system for the unknown coefficients $a_{j}$ is given by  
\begin{equation*}
\left( 
\begin{array}{cccccccc}
1 & 1 & 1 & 1 & 1 & 1 & 1 & \cdot \cdot \cdot  \\ 
1 & -1 & 1 & -1 & 1 & -1 & 1 & \cdot \cdot \cdot  \\ 
&  & 2\cdot 1 & 3\cdot 2 & 4\cdot 3 & 5\cdot 4 & 6\cdot 5 & \cdot \cdot
\cdot  \\ 
&  & 2\cdot 1 & -3\cdot 2 & 4\cdot 3 & -5\cdot 4 & 6\cdot 5 & \cdot \cdot
\cdot  \\ 
&  &  &  & 4! & 5! & 6!/2! & \cdot \cdot \cdot  \\ 
&  &  &  & 4! & -5! & 6!/2! & \cdot \cdot \cdot  \\ 
&  &  &  &  &  & \cdot \cdot \cdot  & \cdot \cdot \cdot 
\end{array}%
\right) .
\end{equation*}%
We see that even in the one-dimensional  case the map $f^{\ast }\rightarrow
\Delta ^{j}f_{|\mathbb{S}^{d-1}}$ is not trivial; this is in fact the map to
the Taylor coefficients 
\begin{equation*}
\left\{ D^{2j}F\left( -1\right) ,D^{2j}F\left( 1\right) \right\} _{j\geq
0}\Longleftrightarrow \left\{ F^{\left( j\right) }\left( 0\right) \right\}
_{j\geq 0}\Longleftrightarrow F\left( z\right) =\sum_{j=0}^{\infty }\frac{%
F^{\left( j\right) }\left( 0\right) }{j!}z^{j}.
\end{equation*}

\subsection{Which polyharmonic functions are in $H^{2}\left( \mathcal{B}%
\right) $ ?}

The polyharmonic function in the one-dimension case satisfy $%
d^{N}P_{N-1}\left( t\right) /dt^{N}=0$ and are polynomials. Their
complexifications $P_{N}\left( z\right) $ belong to the Hardy spaces in the
disc. On the other hand, not all polyharmonic function (of fixed finite
order $N$) belong to the polyharmonic Hardy spaces $H^{2}\left( \mathcal{B}%
\right) .$ In the present section we characterize those polyharmonic
functions which belong to $H^{2}\left( \mathcal{B}\right) .$

We will characterize the polyharmonic functions by means of their boundary
properties. The main point is that the polyharmonic functions in a domain $%
D\subset \mathbb{R}^{d}$ (and more general domains) may be
\textquotedblright parametrized\textquotedblright\ by considering the
Dirichlet problem 
\begin{align*}
\Delta ^{N}u& =0\qquad \text{in }D \\
\Delta ^{j}u& =g_{j}\qquad \text{on }\partial D\quad \text{for }%
j=0,1,2,...,N-1.
\end{align*}%
By means of the classical Green formulas we may find $u\left( x\right) $ for 
$x\in D.$ On the other hand, if the function $u$ which has a  $r-$%
complexification to ball $\mathcal{B}$ of the Klein-Dirac quadric, then by
means of Cauchy type formula (\ref{HuaAronszajnFormula}) we may recover $%
u\left( z,\theta \right) $ using its values $u\left( e^{i\varphi },\theta
\right) $ on the boundary $\partial \mathcal{B}=\mathbb{S}\times \mathbb{S}%
^{d-1}/\mathbb{Z}_{2}.$ Thus we see the boundary values $u\left( e^{i\varphi
},\theta \right) $ provide an alternative parametrization for the
polyharmonic functions. It is essential, as in the one-dimensional case, to
provide a relation between the boundary data $\left\{ g_{j}\right\}
_{j=0}^{N-1}$ and $u\left( e^{i\varphi },\theta \right) .$ 

The following theorem establishes the link between the polyharmonic BVPs and
the boundary values of the elements in $H^{2}\left( \mathcal{B}\right) .$

\begin{theorem}
\label{TpolyharmonicBVP} Let $u$ be a polyharmonic function of order $N\geq 1
$ in the ball $B\subset \mathbb{R}^{d},$ i.e. $\Delta ^{N}u\left( x\right) =0
$ in $B.$ Then the $r-$complexification of $u$ satisfies $u\in H^{2}\left( 
\mathcal{B}\right) $ if and only if 
\begin{equation*}
\Delta ^{j}u_{|\partial B}\in H^{-j}\left( \mathbb{S}^{d-1}\right) ,\qquad 
\text{for }j=0,1,...,N-1,
\end{equation*}%
where $H^{s}\left( \mathbb{S}^{d-1}\right) $ denotes the Sobolev space of
exponent $s.$
\end{theorem}

\proof%
First recall some properties of the polyharmonic operator $\Delta ^{N}.$ The
following operator is defined by 
\begin{equation*}
L_{\left( k\right) }f:=\frac{1}{r^{d+k-1}}\frac{d}{dr}\left[ r^{d+2k-1}\frac{%
d}{dr}\left[ \frac{1}{r^{k}}f\right] \right] ,
\end{equation*}%
see p. $152$ in \cite{okbook}. For $f\left( r\right) =r^{k+2s},$ $k,s\geq 0,$
we obtain 
\begin{align}
L_{\left( k\right) }f\left( r\right) & =\frac{1}{r^{d+k-1}}\frac{d}{dr}\left[
r^{d+2k-1}\frac{d}{dr}\left[ \frac{1}{r^{k}}f\right] \right] =\frac{1}{%
r^{d+k-1}}\frac{d}{dr}\left[ r^{d+2k-1}\frac{d}{dr}\left[ r^{2s}\right] %
\right]   \label{Lkfr} \\
& =4s\left( \frac{d}{2}+k+s-1\right) r^{k+2s-2}.  \notag
\end{align}%
Hence, for every $j\leq s$ holds 
\begin{eqnarray}
L_{\left( k\right) }^{j}\left[ r^{k+2s}\right]  &=&4^{j}s\left( s-1\right)
\cdot \cdot \cdot \left( s-j+1\right) \times   \label{Ljkrk+2s} \\
&&\qquad \qquad \times \left( \frac{d}{2}+k+s-1\right) \cdot \cdot \cdot
\left( \frac{d}{2}+k+s-j\right) r^{k+2s-2j}.  \notag
\end{eqnarray}

Let us put 
\begin{equation}
\gamma _{j,s}^{k}:=L_{\left( k\right) }^{j}\left[ r^{k+2s}\right] _{|r=1}.
\label{gammajs}
\end{equation}%
Obviously, 
\begin{equation*}
\gamma _{j,s}^{k}=0\qquad \text{for }s\leq j-1.
\end{equation*}%
Let the function $u\left( x\right) $ have the expansion 
\begin{equation*}
u\left( x\right) :=\sum_{k=0}^{\infty }\sum_{\ell
=1}^{a_{k}}\sum_{j=0}^{N-1}u_{k,\ell ;j}r^{k+2j}Y_{k,\ell }\left( \theta
\right) ,
\end{equation*}%
and for the boundary data we have the Laplace-Fourier series expansion 
\begin{equation*}
g_{m}\left( \theta \right) :=\sum_{k=0}^{\infty }\sum_{\ell
=1}^{a_{k}}g_{k,\ell }^{m}Y_{k,\ell }\left( \theta \right) .
\end{equation*}%
Hence, by  formula in \cite{okbook}, p. $165$ for the Laplace operator, we
have 
\begin{align*}
g_{m}\left( \theta \right) & =\Delta ^{m}u\left( x\right)
_{|r=1}=\sum_{k=0}^{\infty }\sum_{\ell =1}^{a_{k}}\sum_{j=0}^{N-1}u_{k,\ell
;j}L_{\left( k\right) }^{m}\left[ r^{k+2j}\right] _{|r=1}Y_{k,\ell }\left(
\theta \right)  \\
& =\sum_{k=0}^{\infty }\sum_{\ell =1}^{a_{k}}\sum_{j=0}^{N-1}u_{k,\ell
;j}\gamma _{m,j}^{k}Y_{k,\ell }\left( \theta \right) .
\end{align*}%
By comparing the coefficients we obtain the system of equations 
\begin{equation}
g_{k,\ell }^{m}=\sum_{j=m}^{N-1}u_{k,\ell ;j}\gamma _{m,j}^{k}\qquad \text{%
for }m=0,1,...,N-1.  \label{Systemukl}
\end{equation}%
Hence, the map $u\leftrightarrow \left\{ g_{j}\right\} _{j=0}^{N-1}$ is
given by the infinitely many matrices for $k=0,1,2,...,$ 
\begin{equation*}
U_{k}=\left( 
\begin{array}{ccccc}
1 & 1 & 1 & \cdot \cdot \cdot  & 1 \\ 
0 & \gamma _{1,1}^{k} & \gamma _{1,2}^{k} & \cdot \cdot \cdot  & \gamma
_{1,N-1}^{k} \\ 
0 & 0 & \gamma _{2,2}^{k} & \cdot \cdot \cdot  & \gamma _{2,N-1}^{k} \\ 
\cdot \cdot \cdot  & \cdot \cdot \cdot  & \cdot \cdot \cdot  & \cdot \cdot
\cdot  & \cdot \cdot \cdot  \\ 
0 & 0 & 0 & \cdot \cdot \cdot  & \gamma _{N-1,N-1}^{k}%
\end{array}%
\right) .
\end{equation*}

The result follows from the Lemma \ref{Luklj} below, and by application of
the Cauchy-Bunyakowski-Schwarz inequality to the norm 
\begin{equation*}
\left\Vert u\right\Vert _{H^{2}\left( \mathcal{B}\right)
}^{2}=\sum_{j=0}^{N-1}\sum_{k=0}^{\infty }\sum_{\ell =1}^{a_{k}}\left\vert
u_{k,\ell ;j}\right\vert ^{2}
\end{equation*}%
with equations (\ref{Systemukl}).

\endproof%

We have the following technical lemma.

\begin{lemma}
\label{Luklj}The solution of the system (\ref{Systemukl}) is given by 
\begin{equation*}
u_{k,\ell ;j}=\sum_{p=j}^{N-1}\frac{D_{k}^{p,j}}{D_{k}}\times \frac{%
g_{k,\ell }^{p}}{\left( \frac{d}{2}+k\right) ^{p}}
\end{equation*}%
where $D_{k}=\det U_{k},$ and $D_{k}^{p,j}$ are the appropriate minors of $%
U_{k}$ given by the Cramer's rule. These coefficients satisfy 
\begin{equation*}
\frac{D_{k}^{p,j}}{D_{k}}\rightarrow C_{p,j}\qquad \text{for }k\rightarrow
\infty ,
\end{equation*}%
for appropriate constants $C_{p,j}.$
\end{lemma}

\proof%
In the system (\ref{Systemukl}) we divide the $j-$th row by $\left( \frac{d}{%
2}+k\right) ^{j-1}.$ On the right-hand side we obtain the vector $\left( 
\frac{g_{k,\ell }^{j}}{\left( \frac{d}{2}+k\right) ^{j}}\right) _{j=0}^{N-1},
$ and the coefficients of the system of the left-hand side (\ref{Systemukl})
become equal to 
\begin{equation*}
\frac{\gamma _{j,s}^{k}}{\left( \frac{d}{2}+k\right) ^{j}}\qquad \text{for }%
j\leq s\leq N-1.
\end{equation*}%
From formula (\ref{Ljkrk+2s}) for the constants $\gamma _{j,s}^{k}$ we see
that 
\begin{equation*}
\frac{\gamma _{j,s}^{k}}{\left( \frac{d}{2}+k\right) ^{j}}\rightarrow
F_{j,s}\qquad \text{for }k\rightarrow \infty ,
\end{equation*}%
for some appropriate constants $F_{j,s}$ which are non-zero. This ends the
proof since the minors $D_{p}^{k,j}$ and $D_{k}$ are computed through the
values $\frac{\gamma _{j,s}^{k}}{\left( \frac{d}{2}+k\right) ^{j}}.$

\endproof%

We see that Theorem \ref{TpolyharmonicBVP} provides us with a large class of
functions defined on the ball $B\subset \mathbb{R}^{d}$ which are extendable
to the ball $\mathcal{B}$ of the Klein-Dirac quadric. It is also possible to
consider functions which are in a certain sense polyharmonic of infinite
order, as those studied by Aronszajn, Lelong, Avanissian, and others, cf.
the references in \cite{aron}, \cite{avanissian}, \cite{KounchevRenderPoly}.

\section{Error estimate of Polyharmonic Interpolation and Cubature formulas
\label{Scubature}}

The topic of estimation of quadrature formulas for analytic functions is a
widely studied one. Beyond the classical monographs \cite{krylov}, \cite%
{davisRabinowitz}, we provide further and more recent publications, as \cite%
{bakhvalov}, \cite{gautschi}, \cite{goetz}, \cite{kzaz}, \cite{kowalski}, 
\cite{milova}. No references may be found though for the multivariate case,
for cubature formulas, even in the fundamental monographs as \cite{sobolev}, 
\cite{stroudBook}, \cite{sobolev2}; see also the recent survey \cite{cools}.

Our main framework of Interpolation and Cubature was defined in \cite%
{kounchevRenderArkiv}, \cite{haussmannKounchev}, \cite{kounchevRenderArxiv}.
It has further brought to life the multivariate complexification and the
polyharmonic Hardy spaces.

We will consider first \emph{polyharmonic interpolation,} and as second, 
\emph{polyharmonic Cubature formulas} for approximating integrals of the
type 
\begin{equation*}
\int_{B}f\left( x\right) d\mu \left( x\right) 
\end{equation*}%
over the unit ball $B\subset \mathbb{R}^{d},$ where $\mu \left( x\right) $
is a signed measure of special type. The main purpose of the present section
is to find error estimates for the Interpolation and Cubature in the case of
functions $f\in H^{2}\left( \mathcal{B}\right) .$ Interpolation and
Quadrature are intimately related in the one-dimensional case, and we will
demonstrate similar relation in the multivariate case in our setting.

\subsection{One-dimensional case}

First, following the classical scheme outlined in \cite{davisRabinowitz}, p. 
$303-306,$ (see also chapter $12$ in \cite{krylov}), we will remind how one
finds the error for the classical Quadrature formulas.

Let the points $t_{0},$ $t_{1},$ $...,$ $t_{N}$ belong to the interval $%
\left[ a,b\right] .$ We define 
\begin{equation*}
\omega _{N}\left( z\right) =\left( z-t_{0}\right) \left( z-t_{1}\right)
\cdot \cdot \cdot \left( z-t_{N}\right) .
\end{equation*}%
Let the function $f$ be analytic in a simply connected (open) domain $%
D\subset \mathbb{C}$ containing the interval $\left[ a,b\right] $ with
boundary $\partial D=\Gamma .$ Then the interpolation polynomial $%
P_{N}\left( t\right) $ satisfying $P_{N}\left( t_{j}\right) =f\left(
t_{j}\right) $ for $j=0,1,...,N$ is given by 
\begin{equation*}
P_{N}\left( z\right) =\dsum_{j=0}^{N}f\left( t_{j}\right) \frac{\omega
\left( z\right) }{\omega ^{\prime }\left( t_{j}\right) \left( z-t_{j}\right) 
}=f\left( z\right) -\frac{1}{2\pi i}\dint_{\Gamma }\frac{\omega \left(
z\right) f\left( t\right) }{\omega \left( t\right) \left( t-z\right) }dt,
\end{equation*}%
where $\Gamma $ is considered as a contour oriented counterclockwise. Hence,
the remainder is 
\begin{equation*}
f\left( z\right) -P_{N}\left( z\right) =\frac{1}{2\pi i}\dint_{\Gamma }\frac{%
\omega \left( z\right) f\left( t\right) }{\omega \left( t\right) \left(
t-z\right) }dt.
\end{equation*}%
Now, if $\mu $ is a non-negative Stieltjes measure, say $d\mu \left(
t\right) =w\left( t\right) dt$, the quadrature formula  
\begin{equation*}
\dint_{a}^{b}f\left( t\right) d\mu \left( t\right) \approx
\dsum_{j=0}^{N}\lambda _{j}f\left( t_{j}\right) 
\end{equation*}%
is called\textbf{\ interpolatory quadrature} formulas of degree $N$ if it
satisfies the following equality 
\begin{equation*}
\dsum_{j=0}^{N}\lambda _{j}Q\left( t_{j}\right) =\dint_{a}^{b}P_{N}\left(
t\right) d\mu \left( t\right) 
\end{equation*}%
for every polynomial $Q_{N}$ of degree $\leq N.$ This implies that%
\begin{equation*}
\lambda _{j}=\int_{a}^{b}\frac{\omega \left( t\right) }{\left(
t-t_{j}\right) \omega ^{\prime }\left( t_{j}\right) }d\mu \left( t\right) ;
\end{equation*}%
cf. \cite{davisRabinowitz}, p. $303,$ or Krylov, \cite{krylov}, chapter $12.$
Hence, for the error of such formula we obtain 
\begin{equation*}
E\left( f\right) :=\dint_{a}^{b}\left( f\left( z\right) -P_{N}\left(
z\right) \right) d\mu \left( z\right) =\frac{1}{2\pi i}\dint_{a}^{b}\left(
\dint_{\Gamma }\frac{\omega \left( z\right) f\left( t\right) }{\omega \left(
t\right) \left( t-z\right) }dt\right) d\mu \left( z\right) .
\end{equation*}%
This may be directly estimated by 
\begin{equation*}
\left\vert E\left( f\right) \right\vert \leq \frac{L_{\Gamma }}{2\pi }%
\max_{t\in \Gamma }\left\vert f\left( t\right) \right\vert \times \frac{1}{%
d^{N+1}}\frac{\left( b-a\right) D^{N+1}}{\delta _{\Gamma }}%
\dint_{a}^{b}\left\vert d\mu \left( z\right) \right\vert 
\end{equation*}%
where $d:=\min_{j}\left( \limfunc{dist}\left( t_{j},\Gamma \right) \right) ,$
$D:=\max_{j}\left( \limfunc{dist}\left( t_{j,}\left\{ a,b\right\} \right)
\right) ,$ $\delta _{\Gamma }:=\min \limfunc{dist}\left( \left[ a,b\right]
,\Gamma \right) ,$ and $L_{\Gamma }$ is the length of the contour $\Gamma .$

\subsection{Multivariate Interpolation}

Now we consider the multivariate case.

Let us assume that the number $b$ satisfies 
\begin{equation*}
0<b<1.
\end{equation*}%
Choose the domain $D=\mathcal{B}_{1}$ and a function $f\in H^{2}\left( 
\mathcal{B}\right) ,$ assuming that $f$ has the expansion 
\begin{equation*}
f\left( z,\theta \right) =\sum_{k=0}^{\infty }\sum_{\ell
=1}^{a_{k}}f_{k,\ell }\left( z^{2}\right) z^{k}Y_{k,\ell }\left( \theta
\right) \qquad \text{for }\left\vert z\right\vert <1,\ \text{and }\theta \in 
\mathbb{S}^{d-1}.
\end{equation*}

Let $N\geq 0$ be a fixed integer. We will consider \emph{polyharmonic
interpolation} which has been studied in \cite{haussmannKounchev}. Let the
points $\left\{ r_{k,\ell ;j}\right\} _{j=0}^{N}$ belong to the interval $%
\left[ 0,b\right] .$ We consider the following series: 
\begin{equation}
P_{N}\left( z,\theta \right) =\sum_{k=0}^{\infty }\sum_{\ell
=1}^{a_{k}}p_{k,\ell }\left( z^{2}\right) z^{k}Y_{k,\ell }\left( \theta
\right) ,  \label{PNseries}
\end{equation}%
where and $p_{k,\ell }$ are polynomials of degree $\leq N$ satisfying the
interpolation conditions 
\begin{equation*}
p_{k,\ell }\left( r_{k,\ell ;j}^{2}\right) =f_{k,\ell }\left( r_{k,\ell
;j}^{2}\right) \qquad \text{for }j=0,1,...,N.
\end{equation*}%
We prove below that the series (\ref{PNseries}) is convergent. 

For the remainder of this interpolation we have 
\begin{equation*}
f\left( z,\theta \right) -P_{N}\left( z,\theta \right) =\sum_{k=0}^{\infty
}\sum_{\ell =1}^{a_{k}}\left[ f_{k,\ell }\left( z^{2}\right) -p_{k,\ell
}\left( z^{2}\right) \right] z^{k}Y_{k,\ell }\left( \theta \right) .
\end{equation*}%
Now define as above the functions 
\begin{equation*}
\omega _{k,\ell }\left( z\right) =\left( z-r_{k,\ell ;0}^{2}\right) \left(
z-r_{k,\ell ;1}^{2}\right) \cdot \cdot \cdot \left( z-r_{k,\ell
;N}^{2}\right) ,
\end{equation*}%
and consider the oriented contour $\Gamma \left( t\right) =e^{it}$ for $t\in %
\left[ o,2\pi \right] .$   

For all $z$ with $\left\vert z\right\vert \leq b$ and $\theta \in \mathbb{S}%
^{d-1},$ and obtain the estimate 
\begin{eqnarray}
\left\vert f\left( z,\theta \right) -P_{N}\left( z,\theta \right)
\right\vert  &\leq &\left\vert \sum_{k=0}^{\infty }\sum_{\ell =1}^{a_{k}}
\left[ f_{k,\ell }\left( z^{2}\right) -p_{k,\ell }\left( z^{2}\right) \right]
z^{k}Y_{k,\ell }\left( \theta \right) \right\vert \leq   \label{f-PN} \\
&\leq &\left\vert \sum_{k=0}^{\infty }\sum_{\ell =1}^{a_{k}}\left[
\dint_{\Gamma }\frac{\omega _{k,\ell }\left( z^{2}\right) f_{k,\ell }\left(
\tau ^{2}\right) 2\tau }{\omega _{k,\ell }\left( \tau ^{2}\right) \left(
\tau ^{2}-z^{2}\right) }d\tau \right] z^{k}Y_{k,\ell }\left( \theta \right)
\right\vert   \notag \\
&\leq &\sum_{k=0}^{\infty }\sum_{\ell =1}^{a_{k}}\left\vert \dint_{\Gamma }%
\frac{\omega _{k,\ell }\left( z^{2}\right) f_{k,\ell }\left( \tau
^{2}\right) 2\tau }{\omega _{k,\ell }\left( \tau ^{2}\right) \left( \tau
^{2}-z^{2}\right) }d\tau \right\vert b^{k}\left\vert Y_{k,\ell }\left(
\theta \right) \right\vert   \notag \\
&\leq &\sum_{k=0}^{\infty }\sum_{\ell =1}^{a_{k}}\dint_{\Gamma }\left\vert 
\frac{2^{N+1}f_{k,\ell }\left( \tau ^{2}\right) 2}{\left( 1-b\right) ^{N+2}}%
\right\vert \left\vert d\tau \right\vert \times b^{k}\left\vert Y_{k,\ell
}\left( \theta \right) \right\vert   \notag \\
&\leq &C\frac{2^{N+2}}{\left( 1-b\right) ^{N+2}}\sum_{k=0}^{\infty
}\sum_{\ell =1}^{a_{k}}\dint_{\Gamma }\left\vert f_{k,\ell }\left( \tau
^{2}\right) \right\vert \left\vert d\tau \right\vert \times b^{k}k^{\frac{d-2%
}{2}}  \notag \\
&\leq &C\frac{2^{N+2}}{\left( 1-b\right) ^{N+2}}\sum_{k=0}^{\infty
}\sum_{\ell =1}^{a_{k}}\left\Vert f_{k,\ell }\right\Vert _{H^{2}\left( 
\mathbb{D}\right) }b^{k}k^{\frac{d-2}{2}}.  \notag
\end{eqnarray}%
Since $f\in H^{2}\left( \mathcal{B}\right) $ the last inequality shows,
after application of Cauchy-Bunyakovski-Schwarz inequality, that the series $%
\sum_{k=0}^{\infty }\sum_{\ell =1}^{a_{k}}\left[ f_{k,\ell }\left(
z^{2}\right) -p_{k,\ell }\left( z^{2}\right) \right] z^{k}Y_{k,\ell }\left(
\theta \right) $ is absolutely and uniformly convergent. Hence, the series (%
\ref{PNseries}) representing the function $P_{N}\left( z,\theta \right) $ is
also such.

Above we have outlined the most important arguments for proving the
following:

\begin{theorem}
Let $f\in H^{2}\left( \mathcal{B}\right) .$ Let the points $\left\{
r_{k,\ell ;j}\right\} _{j=0}^{N}$ belong to the interval $\left[ 0,b\right] ,
$ where $b<1.$ Then the function $P_{N}\left( z,\theta \right) $ defined by
the series (\ref{PNseries}) is polyharmonic of order $N+1$ in the ball $%
B_{b}\subset \mathbb{R}^{d}$ and belongs to the polyharmonic Hardy space $%
H^{2}\left( \mathcal{B}_{b}\right) ,$ while the following inequality holds:\ 
\begin{equation*}
\left\Vert P_{N}\right\Vert _{H^{2}\left( \mathcal{B}_{b}\right) }\leq
C_{N,b}\left\Vert f\right\Vert _{H^{2}\left( \mathcal{B}\right) }.
\end{equation*}
\end{theorem}

\subsection{Multivariate Polyharmonic Cubature}

The class of \textbf{pseudo-positive} measures used for our cubature formula 
$C_{N}\left( f\right) $ is now defined in the following way: a signed
measure $\mu $ with support in $B_{R}\subset $ $\mathbb{R}^{d}$ is \emph{%
pseudo-positive with respect to the orthonormal basis} $Y_{k,\ell },\ell
=1,...,a_{k}$, $k\in \mathbb{N}_{0}$ if the inequality 
\begin{equation}
\int_{\mathbb{R}^{d}}h\left( \left\vert x\right\vert \right) Y_{k,\ell
}\left( x\right) d\mu \left( x\right) \geq 0  \label{defpspos}
\end{equation}%
holds for every non-negative continuous function $h:\left[ a,b\right]
\rightarrow \left[ 0,\infty \right) $ and for all $k\in \mathbb{N}_{0}$, $%
\ell =1,2,...,a_{k}.$ Let us note that every signed measure $d\mu $ with
bounded variation may be represented (non-uniquely) as a difference of two
pseudo-positive measures. We refer to \cite{kounchevRenderArkiv} for
instructive examples of pseudo-positive measures.

Let the \emph{pseudo-positive} (signed) measure $d\mu $ be given in the ball 
$B_{b}\subset \mathbb{R}^{d}.$ For all indices $\left( k,\ell \right) $ the
component measures are defined by 
\begin{equation}
d\mu _{k,\ell }\left( r\right) :=\dint_{\mathbb{S}^{d-1}}Y_{k,\ell }\left(
\theta \right) d\mu \left( r\theta \right) \geq 0\qquad \text{for all }r\in %
\left[ 0,b\right] ;  \label{dmukl}
\end{equation}%
here the integral is symbolical with respect to the variables $\theta .$
Rigorously, the component measure $d\mu _{k,\ell }\left( r\right) $ is
defined for the functions $g\left( r\right) $ on the interval $\left[ 0,b%
\right] $ by means of the equality 
\begin{equation*}
\dint_{0}^{b}g\left( r\right) d\mu _{k,\ell }\left( r\right)
:=\dint_{B_{b}}g\left( r\right) Y_{k,\ell }\left( \theta \right) d\mu \left(
x\right) ;
\end{equation*}%
cf. \cite{kounchevRenderArxiv}, \cite{kounchevRenderArkiv}.

In \cite{kounchevRenderArxiv}, \cite{kounchevRenderArkiv}, we have
considered a special type of Cubature formula, the so-called \emph{%
polyharmonic Gauss-Jacobi Cubature formula}. Here however we will consider
more generally, \emph{interpolatory polyharmonic Cubature formulas }and will
prove their convergence and error estimate for them. The case of the annulus
has been considered by us in \cite{kounchevRenderHardyAnnulus}.

Let us fix $\left( k,\ell \right) .$ We assume that there exist points $%
t_{k,\ell ;j},$ $j=0,1,...,N,$ belonging to the interval $\left[ 0,b\right] ,
$ and numbers $\left\{ \lambda _{k,\ell ;j}\right\} _{j=1}^{N},$ such that
the following \textbf{interpolatory quadrature} formula holds: 
\begin{equation}
\dint_{0}^{b}Q\left( r\right) d\mu _{k,\ell }\left( r\right)
=\dsum_{j=0}^{N}\lambda _{k,\ell ;j}Q\left( t_{k,\ell ;j}\right) \qquad 
\text{for every }Q\in V_{k,N};  \label{Quadraturekl}
\end{equation}%
here the set $V_{k,N}$ is given by 
\begin{equation*}
V_{k,N}=\left\{ r^{k+2j}\right\} _{j=0}^{N}.
\end{equation*}%
We define the \textbf{polyharmonic interpolatory Cubature formula} by 
\begin{equation}
C_{N}\left( f\right) :=\sum_{k=0}^{\infty }\sum_{\ell
=1}^{a_{k}}\sum_{j=0}^{N}\lambda _{k,\ell ;j}f_{k,\ell }\left( t_{k,\ell
;j}\right) .  \label{CN}
\end{equation}%
For the interpolation polyharmonic function $P_{N}$ defined in (\ref%
{PNseries}) we obtain equality 
\begin{equation*}
\int_{B}P_{N}\left( x\right) d\mu \left( x\right) =C_{N}\left( P_{N}\right) .
\end{equation*}%
Hence, the remainder of the polyharmonic Cubature formula is given by 
\begin{eqnarray*}
E\left( f\right)  &=&\dint_{B}\left( f\left( r,\theta \right) -P_{N}\left(
r,\theta \right) \right) d\mu \left( z,\theta \right)  \\
&=&\sum_{k=0}^{\infty }\sum_{\ell =1}^{a_{k}}\dint_{0}^{1}\left[ f_{k,\ell
}\left( r^{2}\right) -p_{k,\ell }\left( r^{2}\right) \right] r^{k}d\mu
_{k,\ell }\left( r\right) ,
\end{eqnarray*}%
which implies the estimate 
\begin{align*}
\left\vert E\left( f\right) \right\vert & \leq \sum_{k=0}^{\infty
}\sum_{\ell =1}^{a_{k}}\left\vert \dint_{0}^{1}\left[ f_{k,\ell }\left(
r^{2}\right) -p_{k,\ell }\left( r^{2}\right) \right] r^{k}d\mu _{k,\ell
}\left( r\right) \right\vert  \\
& \leq \sum_{k=0}^{\infty }\sum_{\ell =1}^{a_{k}}\left\vert
\dint_{0}^{1}\dint_{\Gamma }\frac{\omega _{k,\ell }\left( r^{2}\right)
f_{k,\ell }\left( \tau ^{2}\right) 2\tau }{\omega _{k,\ell }\left( \tau
^{2}\right) \left( \tau ^{2}-r^{2}\right) }d\tau \times r^{k}d\mu _{k,\ell
}\left( r\right) \right\vert  \\
& \leq \frac{C_{N}}{\left( 1-b\right) ^{N+2}}\frac{L_{\Gamma }}{2\pi }%
\sum_{k=0}^{\infty }\sum_{\ell =1}^{a_{k}}\left\Vert f_{k,\ell }\right\Vert
_{H^{2}\left( \mathcal{B}\right) }\times \dint_{0}^{1}r^{k}d\mu _{k,\ell
}\left( r\right) .
\end{align*}%
This proves the following result:

\begin{theorem}
Let $f\in H^{2}\left( \mathcal{B}\right) .$ Let the points $\left\{
r_{k,\ell ;j}\right\} _{j=0}^{N}$ belong to the interval $\left[ 0,b\right] $
with $b<1.$ Then the polyharmonic cubature formula defined by (\ref{CN})
with remainder $E\left( f\right) =\int_{B}f\left( x\right) d\mu \left(
x\right) -C_{N}\left( f\right) $ satisfies the following estimate 
\begin{equation*}
\left\vert E\left( f\right) \right\vert \leq \frac{C_{N}}{\left( 1-b\right)
^{N+2}}\frac{L_{\Gamma }}{2\pi }\sum_{k=0}^{\infty }\sum_{\ell
=1}^{a_{k}}\left\Vert f_{k,\ell }\right\Vert _{H^{2}\left( \mathcal{B}%
\right) }\times \dint_{0}^{1}r^{k}d\mu _{k,\ell }\left( r\right) .
\end{equation*}
\end{theorem}

\section{Conclusions}

\begin{enumerate}
\item Our research may be considered as a contribution to the topic of
analytic continuation of solutions to elliptic equations (in particular,
harmonic functions), see the discussion and references to the works of V.
Avanissian, P. Lelong, C. Kiselman, J. Siciak, M. Jarnicki, T. du Cros, on
p. $54-55$ in \cite{aron}, \cite{avanissian}, \cite{Mori98}, \cite{FuMo02}, 
\cite{KoRe08}, \cite{KounchevRenderPoly}. Our construction of $r-$analytic
continuation is applicable to domains as annuli, strips and other domains
with symmetry in $\mathbb{R}^{d}.$ The case of the annulus has been
considered in \cite{kounchevRenderHardyAnnulus}, while the case of strip and
other domains will be considered in \cite{kounchevRenderBook}.

\item The concept of \emph{polyharmonic Hardy spaces} appears to be a new
multivariate concept which differs from the existing approaches in several
complex variables, cf. \cite{stein}, \cite{steinWeiss}, \cite{rudin70}, \cite%
{rudin80}, \cite{coifman}, \cite{sarason1998}, \cite{Shai03}.

\item We have seen that the space of $r-$analytic functions on the
Klein-Dirac quadric provides an useful setting for estimation of the
remainders in Interpolation and Cubature. Although the space of such
functions is $1-1$ mapped to a Hardy space of holomorphic functions of
several complex variables on the Lie ball, our approach has a non-trivial
counterpart on the annulus which is not obtained from $\mathbb{C}^{d}$
constructions, cf. \cite{kounchevRenderHardyAnnulus}. Our approach which is
based on $r-$analytic continuation of solutions to elliptic equations (in
particular, polyharmonic functions) provides non-trivial constructions of
Hardy spaces on complexified annulus, strip and other symmetric domains in $%
\mathbb{R}^{d},$ which are not obtained by the standard approach to
holomorphic functions in $\mathbb{C}^{d}.$
\end{enumerate}

\subsection*{Acknowledgment}

Both authors thank the Alexander von Humboldt Foundation and Grant DO-02-275
with Bulgarian National Foundation.

\end{document}